
\magnification=\magstep1

\baselineskip=1.3\baselineskip

\font\tenmsb=msbm10 \font\sevenmsb=msbm7 \font\fivemsb=msbm5
\newfam\msbfam
\textfont\msbfam=\tenmsb \scriptfont\msbfam=\sevenmsb
\scriptscriptfont\msbfam=\fivemsb
\def\Bbb#1{{\fam\msbfam\relax#1}}

\def\bone{{\bf 1}}

\def\eps{\varepsilon}
\def\R{{\Bbb R}}

\def\P{{\bf P}}

\def\<{\langle}
\def\>{\rangle}

\def\cG{{\cal G}}
\def\cE{{\cal E}}
\def\cV{{\cal V}}
\def\cS{{\cal S}}
\def\cB{{\cal B}}
\def\cT{{\cal T}}
\def\cU{{\cal U}}
\def\cW{{\cal W}}

\def\cL{{\cal L}}
\def\n{{\bf n}}

\def\d{{\bf d}}
\def\v{{\bf v}}
\def\w{{\bf w}}

\def\wt{\widetilde}

\def\prt{{\partial}}
\def\df{{\mathop {\ =\ }\limits^{\rm{df}}}}

\def\ol{\overline}
\def\dist{{\mathop {\rm dist}}}

\def\qed{{\hfill $\square$ \bigskip}}
\def\sqr#1#2{{\vcenter{\vbox{\hrule height.#2pt
        \hbox{\vrule width.#2pt height#1pt \kern#1pt
           \vrule width.#2pt}
        \hrule height.#2pt}}}}
\def\square{\mathchoice\sqr56\sqr56\sqr{2.1}3\sqr{1.5}3}

\input epsf.tex
\newdimen\epsfxsize
\newdimen\epsfysize

\centerline{\bf SHY COUPLINGS}
\footnote{$\empty$} {\rm Research partially supported by
NSF grant DMS-0303310 (KB and ZC).}

\vskip 0.3truein

\centerline{\bf Itai Benjamini{\rm,} Krzysztof
Burdzy {\rm and} Zhen-Qing Chen}

\vskip0.4truein

{\narrower

\noindent
 {\bf Abstract}.
 A pair $(X,Y)$ of Markov processes is called a Markov
 coupling if $X$ and $Y$ have the same transition probabilities
 and $(X,Y)$ is a Markov process.
 We say that a coupling is ``shy'' if
 there exists a (random) $\eps>0$ such that
 $\dist(X_t,Y_t) >\eps$ for all $t\geq 0$.
 We investigate whether shy couplings
 exist for several classes of Markov processes.

  }

\vskip0.5truein

\noindent {\bf 1. Introduction}. The proofs of the main theorems
in two recent papers, [BC2] and [BCJ], contained arguments showing
that certain processes come arbitrarily close to each other, at
least from time to time, as time goes to infinity, with
probability one. The proofs were based on ideas specific to the
models and were rather tedious. We decided to examine several
classes of Markov processes in order to determine the conditions
under which there exists a pair of Markov processes defined on the
same probability space such that each marginal process has the
same transition probabilities and the two processes do not come
close to each other at any time. Although we do not have a
complete solution to this problem, we offer a number of results
whose diversity points to a rich theory. Some of our theorems,
examples and techniques may have interest of their own.

We will focus on two classes of processes---reflected Brownian
motions on Euclidean domains and Brownian motions on graphs. The
second class of processes is really discrete in nature, in the
sense that similar techniques work for random walks on graphs. We
chose these classes of processes because similar processes
appeared in our research in the past.

For a general overview of coupling techniques, see [L].

The rest of the paper is organized as follows. We present basic
definitions and elementary examples in Section 2. Section 3 is
devoted to Brownian motions on graphs. We show that there exists a
shy coupling for Brownian motions on
a graph if all its vertices have degree 3 or higher.
Four examples are also given to illustrate the case when the graph
has some vertices of degree one.
 Section 4 deals with reflected Brownian motions on Euclidean
domains, showing that there exist no shy couplings on $C^1$-smooth
bounded strictly convex domains.

\bigskip
\noindent {\bf 2. Preliminaries and elementary examples}. Unless
specified otherwise, all pairs of processes $(X, Y)$ considered in
this paper will be ``Markov couplings,'' i.e., they will satisfy
the following assumptions.
 \item{(i)} $\{X_t, t\geq 0\}$, $\{Y_t,
t\geq 0\}$ and $\{(X_t,Y_t), t\geq 0\}$ are Markov, and the
transition probabilities for $X$ and $Y$ are identical.
 \item{(ii)} The distribution of $\{X_t, t\geq s\}$ conditional on
$\{(X_s, Y_s) = (x,y)\}$ is the same as the distribution of
$\{X_t, t\geq s\}$ conditional on $\{X_s = x\}$, for all $x,y $
and $s$.

\medskip
Our definition of a Markov coupling is slightly different from
similar concepts in the literature. One could investigate the
question of whether our results hold for ``couplings'' defined in
other ways, for example, whether condition (ii) is essential.
However, we feel that there are more exciting open problems in
this area---see the end of Section 4.

The following elementary discrete-time example shows that there
exist couplings that satisfy (i) but do not satisfy (ii).

\bigskip
\noindent {\bf Example 2.1}. We take $\{0,1\}$ as the state space
of a discrete time Markov process and we let $\{X_k, k\geq 0\}$ be
a sequence of i.i.d. random variables with
$\P(X_k=0)=\P(X_k=1)=1/2$. We will define a process $\{Y_k, k\geq
0\}$ with the same distribution as $\{X_k, k\geq 0\}$. We let
$Y_0$ be independent of $\{X_k, k\geq 0\}$. For $k\geq 1$, we
construct $Y_k$ so that $\P(Y_k=0  \mid X_{k-1} =0)=0.7$ and
$\P(Y_k=1  \mid X_{k-1} =1)=0.7$. Moreover, for every $k\geq 1$, we
make $Y_k$ independent of $X_j$'s for $j< k-1$. It is elementary
to check that $\{X_k, k\geq 0\}$, $\{Y_k, k\geq 0\}$ and
$\{(X_k,Y_k), k\geq 0\}$ are Markov but for $j\geq 1$, the
distribution of $\{Y_k, k\geq j\}$ conditional on $\{Y_{j-1} =0\}$
is not the same as the distribution of $\{Y_k, k\geq j\}$
conditional on $\{(X_{j-1},Y_{j-1}) =(0,0)\}$.

\bigskip

We will assume that the state space $\cS$ for Markov processes $X$
and $Y$ is metric and we will let $\d$ denote the metric. The open
ball with center $x$ and radius $r$ will be denoted $\cB(x,r)$.
The shortest path between two points in $\cS$ will be called a
geodesic. For some pairs of points, there may be more than one
geodesic joining them.

\bigskip
\noindent {\bf Definition 2.2}. A coupling $(X,Y)$ will be called
shy if one can find two distinct points $x$ and $y$ in the state
space with
$$  \P \left( \inf_{0\leq t < \infty}
     \d(X_t, Y_t)>0 \mid X_0=x, Y_0=y \right)>0.
$$

\bigskip

Note that the term ``shy coupling'' is a label for a family of
Markov transition probabilities.

We proceed with completely elementary examples of shy and non-shy
couplings.

\bigskip
\noindent {\bf Examples 2.3}. (i) Let $X$ be a Brownian motion in
$\R^d$ and let $0\ne y\in \R^d$ be a fixed vector. Let $Y_t = X_t
+ y$ for all $t\geq 0$. Then $(X,Y)$ is a shy coupling.

(ii) Let $X$ be a Brownian motion on the unit circle in $\R^2$ and
let $\theta \in (0,2\pi)$ be a fixed number. We define $Y_t$ using
complex notation, $Y_t = e^{i\theta} X_t$ for all $t\geq 0$. Then
$(X,Y)$ is a shy coupling.

(iii) The last two examples can be easily generalized to a wide
class of Markov processes on spaces $\cS$ with a group structure.
If there is a group element $a\ne 0$ such that $a+X$ is a Markov
process having the same transition probabilities as $X$ and
$\inf_{b\in \cS} \d(b, a+b) >0$ then $(X,a+X)$ is a shy coupling.

(iv) Let $X$ and $Y$ be independent Brownian motions in $\R^d$.
Then $(X,Y)$ is a shy coupling if and only if $d\geq 3$.

\medskip
In the sequel, for $a, b \in \R$, $a\wedge b:=\min \{a, \, b\}$,
$a\vee b:=\max\{a, \, b\}$ and $a^+:=a\vee 0$. For $a>0$, $[a]$
denotes the largest integer that does not exceed $a$.

\bigskip
\noindent {\bf 3. Brownian motion on graphs}. In this section, we
will consider processes whose state space is a finite or infinite
graph. More precisely, let $\cG = (\cV,\cE)$ be a graph, where
$\cV$ is the set of vertices and $\cE$ is the set of edges. We
will assume that all vertices have a finite degree, i.e., for
every vertex there are only a finite number of edges emanating
from this vertex, but we do not assume that this number is bounded
over the set of all vertices. We allow an edge to have both
endpoints attached to one vertex. Every vertex will be attached to
at least one edge.
We will identify edges with finite open line
segments (connected subsets of $\R$), with finite and strictly
positive length, and we will identify vertices with topological
endpoints of edges. In this way, we can identify the graph $\cG$
with a metric space $(\cS, \d)$, where $\cS = \cE \cup \cV$, and
$\d(x,y)$ is the shortest path between $x$ and $y$ along the edges
of the graph. We will assume that the length of any edge is
bounded below by $r_0>0$.

Next we will construct ``Brownian motion'' $X$ on $\cS$. See [FW]
for a definition of a general diffusion on a graph. We leave it to
the reader to check that our somewhat informal description of the
process is consistent with the rigorous construction given in
[FW]. By assumption, our process will be strong Markov. Suppose
that $x\in e \in \cE$ and $x$ is not an endpoint of $e$. Recall
that $e$ can be identified with a line segment, say, $e=[0,y]$.
Then $x\in (0,y)$. If $X_0=x$, then the process $X$ evolves just
like the standard one-dimensional Brownian motion until the exit
time from $(0,y)$. Next suppose that $x\in \cS$ is a vertex. Then
there are $n$ edges $e_1, e_2, \dots, e_n$, attached to $x$, with
$n\geq 1$. Choose a small $r>0$ such that the ball $\cB(x,r)$
consists of line segments $I_j$, $j=1,2,\dots, k$, which are
disjoint except that they have one common endpoint $x$. Note that
$n \leq k \leq 2n$, but not necessarily $k=n$, because some edges
may have both endpoints at $x$. We will describe the evolution of
$X$ starting from $x$ until its exit time from $\cB(x,r)$.
Generate a reflected Brownian motion $R$ on $[0,\infty)$, starting
from 0, and kill it at the first exit time from $[0,r]$, denoted
$T_r$. Label its excursions from 0 with numbers $1, 2, \dots, k$,
in such a way that every excursion has a label chosen uniformly
from $\{1, \cdots, k\}$ and independently of all other labels.
Then we define $X_t$ for $t\in[0,T_r]$ so that $\d(X_t, x) = R_t$
and $X_t \in I_j$, where $j$ is the label of the excursion of $R$
from 0 that straddles $t$ (if $R_t=0$ then obviously $X_t =x$).
This defines the process $X_t$ until its exit time from
$\cB(x,r)$. What we said so far and the strong Markov property
uniquely define the distribution of $X$. Note that when the degree
of a vertex $x$ is 1 then $X$ is best described as a process
reflected at $x$. The process $X$ spends zero amount of time at
any vertex.

Recall that we have assumed that the length of all edges is
bounded below by $r_0>0$. Under this assumption, the process
cannot visit an infinite number of vertices in a finite amount of
time. Hence, the above construction defines a process for all
$t\geq 0$. Another consequence of the assumptions that all edges
have length greater then $r_0$ and all vertices have finite degree
is that for any two points in $\cS$ there is only a finite number
of geodesics joining them. It is clear from our construction that
vertices of degree 2 will play no essential role in the paper and
can be ignored. So we will assume without loss of generality that
there are no vertices of degree 2.

\bigskip
\noindent{\bf Theorem 3.1}. {\sl If all vertices of $\cG$ have
degree 3 or higher then there exists a shy coupling for Brownian
motions on $\cS$.}

\bigskip
\noindent{\bf Proof}. We will construct a coupling $(X,Y)$ of
Brownian motions on $\cS$ such that $X$ and $Y$ move in an
independent way when they are far apart and they move in a
``synchronous'' way when they are close together. Clearly,
independent processes do not form a shy coupling on a finite
graph. Remark 3.2 below explains why it is hard, perhaps
impossible, to construct a ``synchronous'' shy coupling.

For any $x,y\in \cS$ with $\d(x,y) > r_0/4$, we will define
$(X_t,Y_t)$ starting from $(X_0,Y_0)=(x,y)$, for $t\in[0,\tau]$,
where $\tau$ is a random time depending on $x$ and $y$. Then we
will explain how one can define $(X_t,Y_t)$ for $t\in[0,\infty)$
by pasting together different pieces of the trajectory.

(i) Recall that the length of any edge is at least $r_0>0$. First
suppose that $\d(x,y) \geq 3r_0/4$. Then we let $\{(X_t,Y_t),
t\in[0,\tau]\}$ be two independent copies of Brownian motion on
$\cS$ and we let $\tau=\inf\{t>0: \d(X_t,Y_t) =r_0/2\}$.

(ii) Next suppose that $x,y\in \cS$ are such that $\d(x,y)\in(
r_0/4, 3r_0/4)$, and none of these points is a vertex. Let
$$ \sigma (r)= {(4|r| -r_0 )^+ \over r_0} \wedge 1,
 \eqno(3.1)
$$
and $B$ and $B'$ be independent Brownian motions on $\R$ starting
from the origin. Let $U_t=B_t$ and
 $$dV_t = \sqrt{ 1- \sigma^2 (U_t - V_t)} dB_t +
 \sigma (U_t - V_t) dB'_t , \eqno(3.2)
$$
with $V_0=v_0=\d(x,y)> r_0/4$.
Then, if we write $Z_t = V_t - U_t$, we obtain
 $$Z_t = v_0 + \int_0^t (\sqrt{ 1- \sigma^2 (Z_s)}-1) dB_s
 + \int_0^t \sigma(Z_s) dB'_s,$$
and for $Z'_t \df Z_t - r_0/4$,
 $$Z'_t = v_0 -r_0/4
 + \int_0^t (\sqrt{ 1- \sigma^2 (Z'_s+r_0/4)}-1) dB_s
 + \int_0^t \sigma(Z'_s+r_0/4) dB'_s.$$
So
$$ Z_t'= v_0 -r_0/4 +\int_0^t \gamma (Z_s') dW_s,
$$
where
 $$ \gamma(r) \df (\sqrt{ 1- \sigma^2 (r+r_0/4)}-1)^2
 + \sigma^2(r+r_0/4)
$$
and $W$ is a Brownian motion on $\R$ with $W_0=0$.
The process $Z'$ has the same distribution as
$$ t\mapsto  v_0 -r_0/4 +W_{\tau_t},
$$
where
$$ \tau_t:= \inf\left\{ s>0: \int_0^s \gamma( v_0 -r_0/4 + W_s)^{-2} ds>t \right\}.
$$
Note that for small $r>0$, $\gamma (r)=O(r^2)$ and so in particular
$\int_{0+} \gamma(r)^{-2} dr = \infty$.
Thus by Lemma V.5.2 of [KS],
$$ \int_0^{T_0} \gamma (v_0 -r_0/4 + W_s)^{-2} ds =\infty
$$
almost surely, where $T_0 =\inf\{t>0: v_0 -r_0/4 + W_s=0\}$.
We conclude that $Z'$ never hits 0; in other words, $Z$ never reaches
$r_0/4$.

Suppose that $X_0=x$, $Y_0=y$ and recall that we have assumed that
$\d(x,y)\in( r_0/4, 3r_0/4)$. Suppose that $x\in e_1 \in \cE$ and
$y\in e_2 \in \cE$. Let $e_1 \setminus \{x\}$ consist of two line
segments $e_1^\ell$ and $e_1^r$, with $e_1^r$ being the one closer
to $y$. Similarly, $e_2 \setminus \{y\}$ consists of two line
segments $e_2^\ell$ and $e_2^r$, and $e_2^\ell$ is closer to $x$.
We will define $X$ and $Y$ on an interval $[0,\tau]$ to be
specified later.
We define $X_t$ on $e_1$ to be such that $\d(X_t, x) = |U_t|$ and
$X_t \in e_1^\ell$ if and only if $U_t <0$. We define the process
$Y_t$ on $e_2$ by  conditions $\d(Y_t, y) = |V_t-v_0|$ and  $Y_t
\in e_2^\ell$ if and only if $V_t <v_0$.
 We let $\tau$ be the first time $t>0$ that $X_t$ or $Y_t$ is at a
vertex, or $\d(X_t,Y_t) = 3r_0/4$. We see that over the interval
$[0, \tau )$, the distance between $X_t$ and $Y_t$ remains in the
interval $(r_0/4, \, 3r_0/4)$.

\medskip

(iii) This part of our argument is based on the ``skew Brownian
motion.'' The skew Brownian motion $U$ is a real-valued diffusion
which satisfies the stochastic differential equation
 $$ U_t =  B_t + \beta L^U_t, \eqno(3.3) $$
where $B$ is a given Brownian motion with $B_0=0$, $\beta \in[-1,1]$
is a fixed constant and $L^U$ is the symmetric local time of $U$
at $0$, i.e.,
 $$ L^U_t = \lim_{\eps \to 0} {1\over 2\eps} \int_0^t
 \bone_{(-\eps, \eps)} (U_s ) ds
 \, . \eqno(3.4)
 $$
The existence and uniqueness of a strong solution to (3.3)-(3.4)
was proved in [HS]. In the special case of $\beta =1$, the
solution to (3.3) is the reflected Brownian motion. An alternative
way to define the skew Brownian motion is the following. Consider
the case $\beta > 0$. Take a standard Brownian motion $B_t'$ and
flip every excursion of $B_t'$ from $0$ to the positive side with
probability $\beta$, independent of what happens to other
excursions (if an excursion is on the positive side, it remains
unchanged). The resulting process has the same distribution as $U$
defined by (3.3)-(3.4). For more information and references, see
recent papers on skew Brownian motion, [BC1] and [BK].

Suppose that $x,y\in \cS$, $\d(x,y)\in(r_0/4, 3r_0/4)$, and $x$ is
a vertex. Note that $y$ is not a vertex. By assumption, the degree
$k$ of vertex $x$ is 3 or greater. Suppose that $B_t$ is a
Brownian motion on $\R$ and let $U$ be a solution to (3.3)-(3.4),
with $\beta$ defined by $(1-\beta)/(1+\beta) = k-1$. Note that
$\beta<0$. We label negative excursions of $U$ from 0 with numbers
$1, 2, \dots, k-1$, in such a way that every excursion has a label
chosen uniformly from this set and independently of all other
labels.

Suppose that $B'$ is a Brownian motion independent of $B$. Recall
the definition of the function $\sigma$ and the process $V$ given
in (3.1) and (3.2), respectively, with $V_0=v_0=\d(x,y)$. Then,
if we write $Z_t = V_t - U_t$, we obtain
 $$Z_t = v_0 -\beta L^U_t
 +\int_0^t (\sqrt{ 1- \sigma^2 (Z_s)}-1) dB_s
 + \int_0^t \sigma(Z_s) dB'_s,$$
and for $Z'_t \df Z_t - r_0/4$,
 $$Z'_t = v_0 -r_0/4 -\beta L^U_t
 + \int_0^t (\sqrt{ 1- \sigma^2 (Z'_s+r_0/4)}-1) dB_s
 + \int_0^t \sigma(Z'_s+r_0/4) dB'_s.$$
We have already pointed out that for small $r>0$,
 $$ \gamma(r) \df (\sqrt{ 1- \sigma^2 (r+r_0/4)}-1)^2
 + \sigma^2(r+r_0/4) =O(r^2).
$$
Since $\beta<0$, the process $- \beta L^U_t$ is nondecreasing.
These observations and the argument used in the first half of (ii)
 imply that $Z'$ never hits 0, i.e., $Z$ never reaches $r_0/4$.

The ball $\cB(x,3r_0/4)$ consists of line segments $I_j$,
$j=1,2,\dots,k$. We assume that $I_k$ is the line segment
containing $y$. We define $X_t$ on these line segments so that
$\d(X_t, x) = |U_t|$. If $U_t >0$ then $X_t \in I_k$. If $U_t < 0$
then $X_t \in I_j$, where $j$ is the label of the excursion of $U$
straddling $t$. Suppose that $y\in e\in \cE$. Let $e \setminus
\{y\}$ consist of two line segments $e^\ell$ and $e^r$, with
$e^\ell$ being the one closer to $x$. We define $Y_t$ on $e$ by
$\d(Y_t,y)=|V_t-v_0|$.
We let $Y_t \in e^\ell$ if and only if $V_t <v_0$. We let $\tau$ be the
infimum of $t$ such that $Y_t$ is at a vertex, or $\d(X_t,Y_t) =
3r_0/4$. Observe that over the interval $[0, \tau )$, the distance
between $X_t$ and $Y_t$ remains in the interval $(r_0/4, \,
3r_0/4)$.

\medskip

Now we will define the process $(X_t,Y_t)$ for all $t\geq 0$,
assuming that $X_0 =x$, $Y_0=y$ and $\d(x,y) > r_0/4$. We use one
of the parts (i)-(iii) of the proof to define the process $(X,Y)$
on an interval $[0,\tau_1]$. Then we proceed by induction. Suppose
that the process has been defined on an interval $[0,\tau_k]$ and
$\d(X_{\tau_k},Y_{\tau_k})>r_0/4$. Then we use the appropriate
part (i)-(iii) of the proof to extend the process, using the
strong Markov property at $\tau_k$, to an interval
$[0,\tau_{k+1}]$. It is easy to see that $\tau_k\to\infty$ a.s.,
so the process $(X_t,Y_t)$ is defined for all $t\geq 0$. It is
straightforward to check that $\{X_t, t\geq 0\}$ and $\{Y_t, t\geq
0\}$ are Brownian motions on $\cS$ and $(X,Y)$ is a shy Markov
coupling, as defined in Section 2.
 \qed

\bigskip
\noindent{\bf Remark 3.2}. One may wonder whether it is possible
to construct the shy coupling in the proof of Theorem 3.1 using
the skew Brownian motion in such a way that the distance between
$X$ and $Y$ does not change on time intervals where both processes
stay away from vertices (``synchronous coupling''). This idea
works well for many graphs but runs into technical problems when
we have a configuration similar to that in
Fig.~3.1, with many
geodesics joining two vertices. Suppose that the lengths of edges in
Fig.~3.1
are chosen so that there are 6 geodesics between $x$
and $y$. We will argue that if $X_0=x$ and $Y_0=y$ then for small
$t>0$ the distance between $X$ and $Y$ has to decrease. This is
because $X$ will have to move towards $z$ with probability $5/7$,
and $Y$ will have to move towards $z$ with probability $1/3$.
Since $5/7 + 1/3>1$, $X$ and $Y$ will find themselves on a
geodesic from $x$ to $y$, moving away from their starting points
towards each other, with positive probability. For this reason, we
could not find a ``synchronous'' shy coupling based on the skew
Brownian motion.

\bigskip \vbox{ \epsfxsize=3.0in
  \centerline{\epsffile{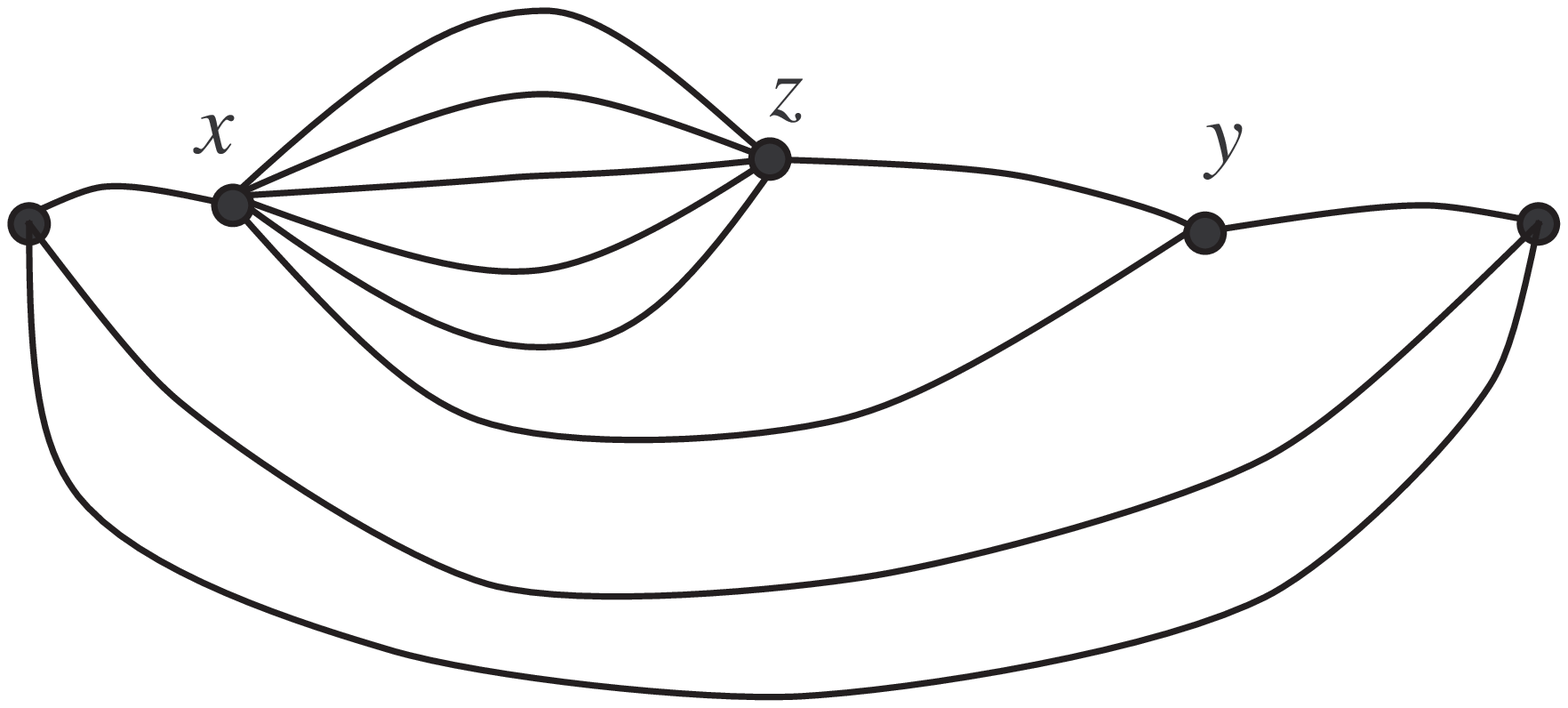}}

\centerline{Figure 3.1.} }
\bigskip

The rest of this section is devoted to graphs that have at least
one vertex of degree 1 (recall that vertices of degree 2 can be
ignored and we assume that $\cG$ does not contain any of them). We
do not have a general theorem covering all graphs with some
vertices of degree 1 but we have four examples illustrating some
special cases.

\bigskip
\noindent{\bf Example 3.3}. This example is similar to Examples
2.3 (i)-(iii). Suppose that there exists an isometry $I: \cS\to
\cS$ such that $\inf_{x\in\cS} \d(x, I(x))>0$. It is not hard to
show that this holds if the isometry has no fixed points, i.e., if
there does not exist $x\in \cS$ with $I(x) =x$. If such an
isometry exists, then we can first construct the process $X$ and
then take $Y_t = I(X_t)$ for all $t\geq 0$. Obviously, thus
constructed coupling $(X,Y)$ is shy. Fig.~3.2 shows that a graph
with some vertices of degree 1 may have this property.

\bigskip \vbox{ \epsfxsize=1.5in
  \centerline{\epsffile{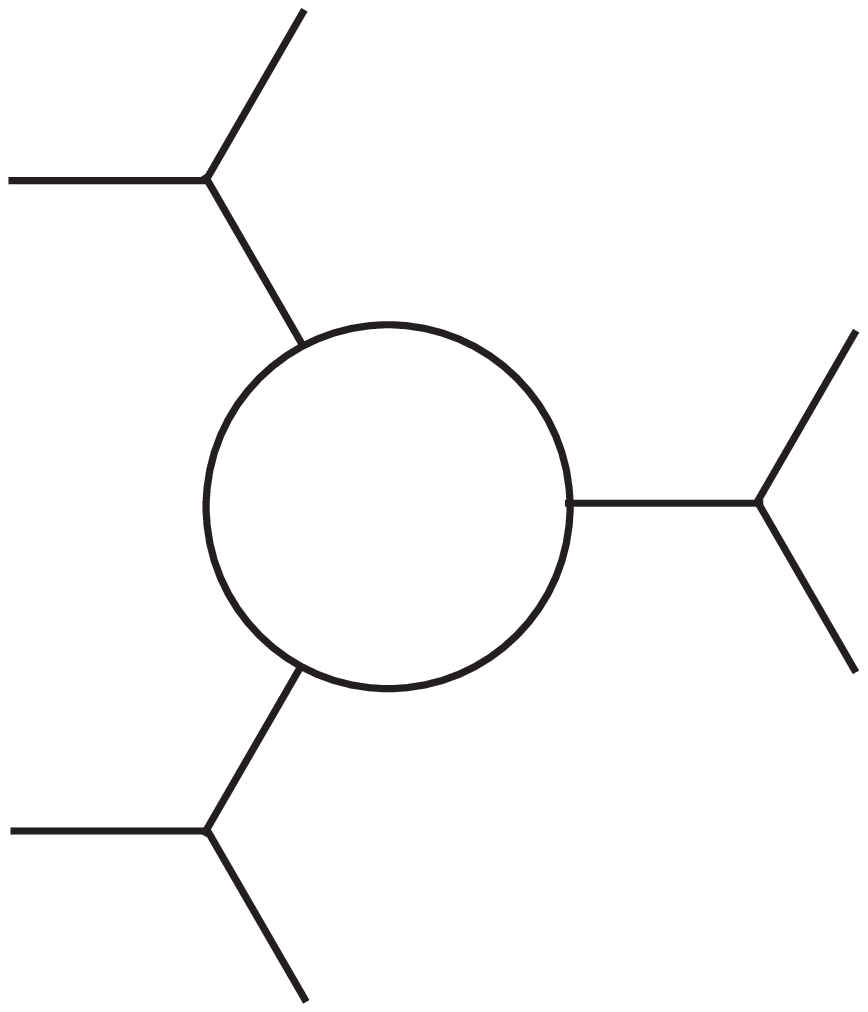}}

\centerline{Figure 3.2.} }
\bigskip

Our next lemma is a large deviations-type estimate. Recall that
all edges are at least $r_0>0$ units long, by assumption.

\bigskip

Let $X$ be Brownian motion on a graph $\cS$. For $A\subset \cS$, define
$T_A = \inf\{t\geq 0: X_t\in A\}$.

\bigskip

\noindent{\bf Lemma 3.4}. {\sl Assume that the degrees of all
vertices are bounded above by $m_0$. There exist constants
$c_0>0$, $t_0<\infty$ depending on $r_0$ only such that for
$t\in(0,t_0)$ and $r>0$ with $r^2> t$ and $\cB(x,r)^c \ne
\emptyset$, we have
 $$
 \left( {c_0\over m_0} \sqrt{t \over 2 \pi}  {1\over 2 r}
 \right)^{[r/r_0]} \exp\left(-{r^2\over 2t}\right)
 \leq \P(T_{\cB(x,r)^c} < t) \leq
 (m_0^{[r/r_0]})!
 \sqrt{2t \over \pi}
 {1\over r} \exp\left(-{r^2\over 2t}\right).$$

}

\bigskip

In applications of the above estimates, $r$ will be ``fixed'' and
then a $t$ (much smaller than $r^2$) will be chosen. For this
reason, we did not try to optimize the non-exponential
factors---most likely they are not best possible.

\bigskip
\noindent{\bf Proof}. (i) First we will prove the lower bound. We
start with some preliminary estimates.

Suppose that $B$ is a Brownian motion on $\R$ with $B_0=0$ and let
$T^B_r=\inf\{t\geq 0: B_t =r\}$. Then $\P(T^B_r < t) = 2 {1\over
\sqrt{2 \pi t}} \int_r^\infty e^{-u^2/2t} du$ for $r>0$. The
following inequalities are well known (see Problem 9.22 on page
112 of [KS]):
 $${r \over 1 + r^2} e^{-r^2/2} \leq
 \int_r^\infty e^{-u^2/2} du \leq {1 \over r} e^{-r^2/2}.
 $$
So for $r\geq 1$,
 $$ {1 \over 2 r} e^{-r^2/2} \leq
 \int_r^\infty e^{-u^2/2} du \leq {1 \over r} e^{-r^2/2}.
 $$
By scaling we obtain for $t\leq r^2$,
 $$ \sqrt{t \over 2 \pi}
 {1\over r} \exp\left(-{r^2\over 2t}\right)
 \leq \P (T^B_r<t)= \sqrt{{2\over \pi}} \int_{r/\sqrt{t}}^\infty
 e^{-v^2/2} dv
 \leq \sqrt{2t \over \pi}
 {1\over r} \exp\left(-{r^2\over 2t}\right).\eqno(3.5)$$
For $r\leq 1$, we have a trivial lower bound $ \int_1^\infty
e^{-u^2/2} du \df c_0>0$. For $t\geq r^2$, the same upper bound
holds but the lower bound has to be replaced with a trivial bound
$c_1=\sqrt{2\over \pi} c_0>0$.

Let $t_0< \infty$ be the largest real such that
$$ 1 - 2 \exp(- r_0^2/ (2t_0))\geq 1/2 .
$$
We will derive an estimate for $\P(T^B_r < t\land T^B_{-r_0/2})$,
for $r\geq r_0/2$ and $t< t_0$. If $r^2<t_0$, then $r\leq c_2 r_0$
for some constant $c_2 < \infty$. In this case, it follows easily
from the support theorem for Brownian motion that $\P(T^B_r <
t\land T^B_{-r_0/2})> c_3>0$ for every $t\in (r^2, \, t_0)$.

Now suppose that $t\leq r^2\land t_0$. If Brownian motion hits
$-r_0/2$ and then it reaches $r$ in $t$ seconds or less, it has to
go from level $-r_0/2$ to level $r$ in $t$ seconds or less. Hence,
by the strong Markov property applied at $T^B_{-r_0/2}$,
 $$\eqalignno{
 \P&(T^B_r < t\land T^B_{-r_0/2}) \geq
 \P(T^B_r < t) - \P(T^B_{r+r_0/2} < t)\cr
 &\geq \sqrt{t \over 2\pi}
 {1\over r} \exp\left(-{r^2\over 2t}\right)
 - \sqrt{2t \over \pi}
 {1\over r+r_0/2} \exp\left(-{(r+r_0/2)^2\over 2t}\right)\cr
 &\geq \sqrt{t \over 2 \pi}
 {1\over r} \exp\left(-{r^2\over 2t}\right)
 - 2\sqrt{t \over 2 \pi}
 {1\over r} \exp\left(-{r^2\over 2t}\right)
 \exp\left(-{ r_0^2\over 2t}\right)\cr
 & = \sqrt{t \over 2 \pi}  {1\over r}
 \exp\left(-{r^2\over 2t}\right)
 \left[ 1 - 2 \exp\left(-{ r_0^2\over 2t}\right)\right]\cr
 &\geq \sqrt{t \over 2 \pi}  {1\over 2 r}
 \exp\left(-{r^2\over 2t}\right).&(3.6)
 }$$

\medskip

Suppose that $y_1$ is a vertex, $\d(y_1,y_2) = r_0/2$, and $y_3$
lies between $y_1$ and $y_2$ so that $\d(y_1,y_3) + \d(y_3, y_2) =
r_0/2$. Let $\d(y_1,y_3) =r_1 \in[0,r_0/2]$ and suppose that $X_0
= y_3$. The process $R_t \df \d(X_t, y_1)$ is a one-dimensional
reflected Brownian motion with $R_0=r_1$, at least until it
reaches $r_0$. Let $\{z_1, z_2, \dots, z_k\}$ be the set of all
points with $\d(z_j, y_1) = r_0/2$ ($y_2$ is one of these points).
Let $T^R_{r_0/2} = \inf\{t\geq 0: R_t = r_0/2\}$. If $X$ visited
$y_1$ before $T^R_{r_0/2}$ then it is at $y_2$ at time
$T^R_{r_0/2}$ with probability $1/k$, by symmetry. Since $k\leq
m_0$, the probability that $X$ starts at $y_3$ and reaches $y_2$
in $s$ seconds or less is greater than or equal to the probability
that reflected Brownian motion that starts from $r_1$ reaches
$r_0/2$ in $s$ seconds or less, divided by $m_0$. The last
probability is bounded below by the analogous probability for the
non-reflected Brownian motion, so using (3.5) we obtain,
 $$\P(T^X_{y_2} \leq s \mid X_0 = y_3) \geq
 {1\over m_0} \sqrt{s \over 2\pi}
 {1 \over r_0/2-r_1}
 \exp\left(-{(r_0/2 - r_1)^2\over 2s}\right), \eqno(3.7)
 $$
for $s \leq (r_0/2 - r_1)^2$. If $s \geq (r_0/2 - r_1)^2$, the
bound is $c_1/m_0$. Note that these estimates hold for all $r_1
\in[0,r_0/2]$, including $r_1 =0$.

Consider any $x'\in \prt \cB(x,r)$ and let $\Gamma\subset \cS$ be
a geodesic connecting $x$ and $x'$. Suppose that $\Gamma $
contains some vertices and denote them $x_1, x_2, \dots, x_k$, in
order in which they lie on $\Gamma$, going from $x$ to $x'$. Let
$x_0 = x$ and $x_{k+1} = x'$. If there is a vertex closer to $x_0$
than $r_0/2$ and it is not $x_1$ then we let $y_0$ be the point at
the distance $r_0/2$ from that vertex, between $x_0$ and $x_1$.
For every $x_j$, $j\geq 1$, we let $y_j\in \Gamma$ be the point
$r_0/2$ away from $x_j$, between $x_j$ and $x_{j+1}$.

Let $z_j$, $j=1, \dots, m_1$,  be the sequence of all points $x_j$
and $y_j$, in the order in which they appear on $\Gamma$ from $x$
to $x'$, including $x$ and $x'$. Note that $\sum_{1\leq j \leq
m_1-1}\d(z_j, z_{j+1})=r$. By the strong Markov property applied
at the hitting times of $z_j$'s, the probability that $X$ starting
from $x$ will hit $x'$ in $t$ seconds or less is bounded below by
$\prod_{j=1}^{m_1-1} p_j$, where $p_j$ is the probability that $X$
starting from $z_j$ will hit $z_{j+1}$ in $t_j$ seconds or less,
and $t_j = t \d(z_j, z_{j+1})/r$.

If $z_j$ is a vertex or $x_0$ then, by (3.7), if $t_j \leq \d(z_j,
z_{j+1})^2$,
 $$p_j \geq
 {1\over m_0} \sqrt{t_j \over 2\pi}
 {1 \over \d(z_j, z_{j+1})}
 \exp\left(-{\d(z_j, z_{j+1})^2\over 2t_j}\right).\eqno(3.8)$$
If $t_j \geq \d(z_j, z_{j+1})^2$ then
 $$p_j \geq c_1/m_0 \geq (c_1/m_0)
 \exp\left(-{\d(z_j, z_{j+1})^2\over 2t_j}\right).\eqno(3.9)$$
For other $z_j$'s we use (3.6) to obtain, for $t_j \leq \d(z_j,
z_{j+1})^2$,
 $$p_j \geq
 \sqrt{t_j \over 2\pi}
 {1 \over 2\d(z_j, z_{j+1})}
 \exp\left(-{\d(z_j, z_{j+1})^2\over 2t_j}\right),\eqno(3.10)$$
and for $t_j \geq \d(z_j, z_{j+1})^2$,
  $$p_j \geq c_3 \geq c_3
 \exp\left(-{\d(z_j, z_{j+1})^2\over 2t_j}\right).\eqno(3.11)$$

The product of exponential factors on the right hand sides of
(3.8)-(3.11) is equal to
 $$ \prod_{j=1}^{m_1-1} \exp \left( - {\d(z_j,z_{j+1})^2
 \over 2 t_j}\right )
 =  \exp \left( - {r^2 \over 2 t}\right ).\eqno(3.12)
 $$

If $t_j \geq \d(z_j,z_{j+1})^2 $ then the non-exponential factor
in (3.9) is $c_1/m_0$ and it is $c_3$ in (3.11). The
non-exponential factors in (3.8) and (3.10) are bounded below by
 $${1\over m_0} \sqrt{t_j \over 2 \pi}  {1\over 2 \d(z_j,z_{j+1})}
 = {1\over m_0} \sqrt{t \d(z_j, z_{j+1})/r \over 2 \pi}
 {1\over 2 \d(z_j,z_{j+1})}
 \geq {1\over m_0} \sqrt{t \over 2 \pi}
 {1\over 2 r}.$$
We conclude that the product of non-exponential factors in
(3.8)-(3.11) is bounded below by
 $$\left( {c_4\over m_0} \sqrt{t \over 2 \pi}  {1\over 2 r}
 \right)^{m_1-1}
 \geq \left( {c_4\over m_0} \sqrt{t \over 2 \pi}  {1\over 2 r}
 \right)^{[r/r_0]}.
 $$
This
combined
 with (3.12) gives the lower bound in the lemma.

(ii) Next we will prove the upper bound. Let $\{\Gamma_j\}$ be the
family of all Jordan arcs in $\cS$ linking $x$ with $\prt
\cB(x,r)$. The number of edges in $\cB(x,r)$ is bounded by $m_2 =
m_0^ {[r/r_0]}$ so the number of $\Gamma_j$'s is bounded by $m_3=
m_2!$. The length of any $\Gamma_j$ is at least $r$.

Consider some $\Gamma_k$. We will define a process $R^k_t$ that
measures the distance from $X_t$ to $x$ along $\Gamma_k$, in a
sense. We will ``erase'' excursions away from $\Gamma_k$ and loops
as follows. For $t>0$, let $\ell(t) = \sup\{s\leq t: X_s \in
\Gamma_k\}$. Let $\cV_k= \cV\setminus \Gamma_k$, i.e., $\cV_k$ is
the set of vertices that do not belong to $\Gamma_k$. Let $T_1=0$,
 $$\eqalign{
 S_j &= \inf\{t\geq T_j: X_t \in \cV_k\}, \quad j \geq 1,\cr
 T_{j+1} & = \inf\{t\geq S_{j}: X_t = X_{\ell(S_{j})}\}, \quad j\geq 1.
 }$$
If $t\in [S_j, T_{j+1}]$ for some $j\geq 1$, we let $R^k_t$ be the
distance from
$X_{\ell(S_{j})}$
 to $x$ along $ \Gamma_k$. For other $t$, we let $R^k_t$ be the
distance from $X_{\ell(t)}$ to $x$ along $\Gamma_k$. The process
$R^k_t$ is a time-change of reflected Brownian motion, that is, it
is reflected Brownian motion ``frozen'' on time intervals when $X$
is outside $\Gamma_k$ (and some other intervals). Hence, the
probability that $R^k_t$ reaches $r$ in $t$ seconds or less is
less than the right hand side of (3.5). Note that one of the
processes $R^k_t$ must be at the level $r$ at the time when $X$
hits $\prt \cB(x,r)$. Hence, an upper bound on the probability in
the statement of the lemma is the product of the right hand side
of (3.5) and $m_3$.
 \qed

\bigskip
\noindent{\bf Example 3.5}. Suppose that the graph $\cS$ is
compact and has the following structure. For some $x\in \cS$, the
set $\cS\setminus \{x\}$ is disconnected and consists of a finite
number of disjoint finite trees $\cT_1, \cT_2, \dots, \cT_k$, and
a graph $\cU$ (not necessarily a tree). We say that a vertex of a
tree is a leaf if it has degree 1. Assume that for some
$r_1>r_2>0$ and every leaf $y\ne x$ of any tree $\cT_j$ we have
$\d(x,y) \geq r_1$, and for any $z\in \cU$, $\d(x,z) \leq r_2$
(see, for example, Fig.~3.3). Suppose that $(X,Y)$ is a coupling
of Brownian motions on $\cS$. We will show that
$(X,Y)$ is not a shy coupling.

\bigskip \vbox{ \epsfxsize=2.5in
  \centerline{\epsffile{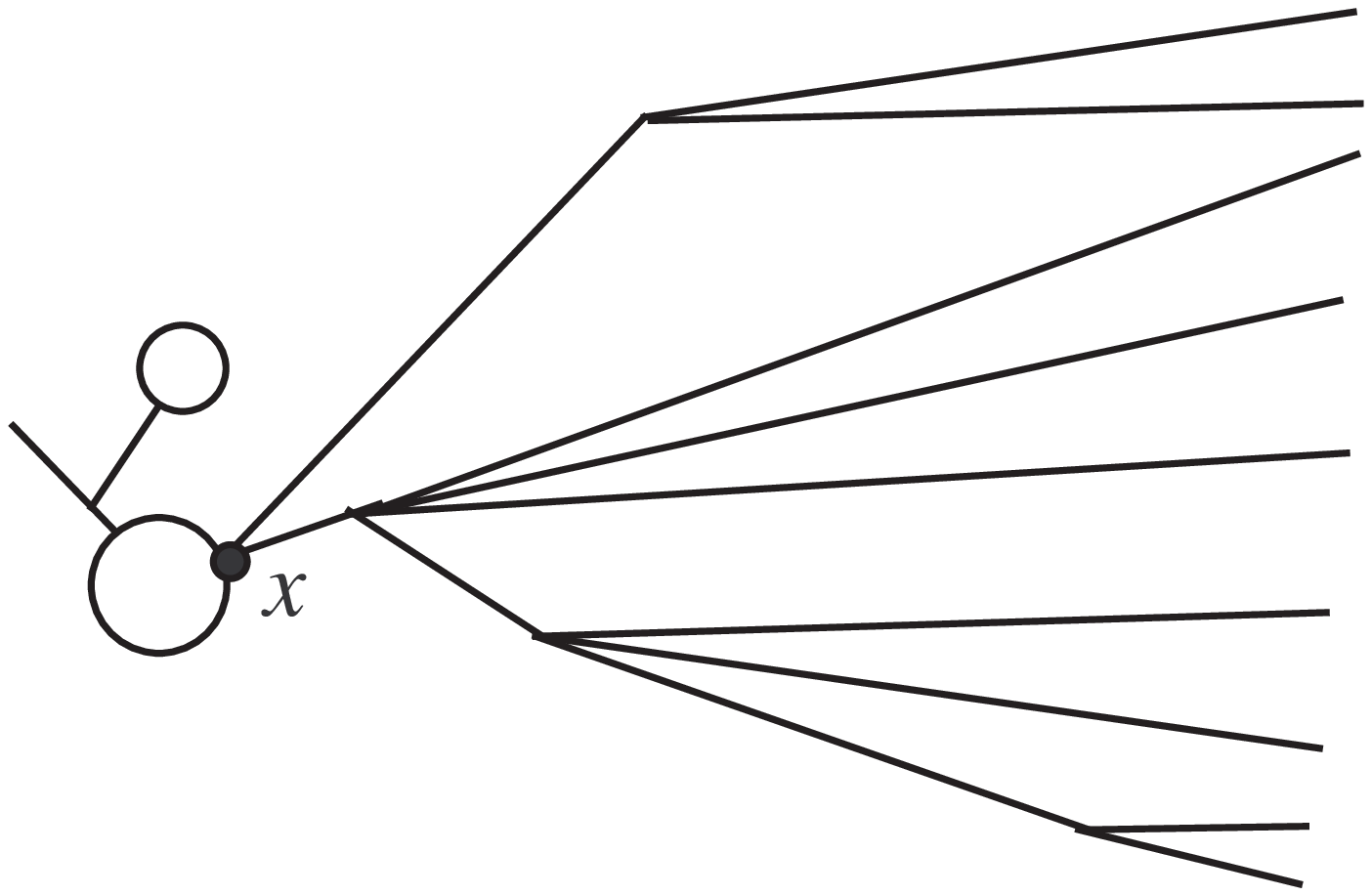}}

\centerline{Figure 3.3.} }
\bigskip

(i) Let $\wt \cT = \bigcup _j \cT_j \setminus \{x\}$. First, we
will show that there exist $p_1>0$ and a stopping time
$T_1<\infty$ such that with probability greater than $p_1$, either
$X_{T_1}=x$ and $Y_{T_1}\in \wt \cT$, or $Y_{T_1}=x$ and
$X_{T_1}\in \wt \cT$. Let $r_3 \in (r_2, r_1)$ and $\cW =
\{y\in\wt \cT: \d(y,x)\geq r_3\}$. It is easy to see that $X$ is
recurrent so $T_2 \df \inf\{t\geq 0: X_t \in \cW\} < \infty$ a.s.
Suppose first that $Y_{T_2} \in \wt \cT\cup \{x\}$ and let $T_3
=\inf\{t\geq T_2: X_t = x \hbox {  or } Y_t =x\}$. Then $T_1 =
T_3$ has the properties stated above.

Next suppose that $Y_{T_2} \in \cU$ and let $T_4 = \inf\{t\geq T_2:
X_t = x\}$ and $T_5 = \inf\{t\geq T_2: Y_t = x\}$. By Lemma 3.4,
for some $p_2, t_0>0$ and all $y\in \cU$ and $z\in \cW$,
$$
  \P(T_5< T_2 + t_0\mid X_{T_2}=z, Y_{T_2}=y)
   > \P(T_4 < T_2 + t_0\mid X_{T_2}=z, Y_{T_2}=y) + p_2.
$$
Hence, $\P(T_5 < T_4\mid X_{T_2}=z, Y_{T_2}=y) > p_2$, and it follows
that we can take $T_1 = T_5$ on the event $\{T_5 < T_4\}$. This
completes the proof of our claim, with $p_1=p_2$.

(ii)
Recall that the length of any edge is bounded below by $r_0>0$.
Fix an arbitrarily small $\eps\in(0, r_0/3)$. We will show in the
remaining part of the proof that $X$ and $Y$ come within $\eps$
distance to each other in finite time almost surely, which will
then imply that $(X,Y)$ is not a shy coupling.

The rest of the proof is based on an inductive argument. We will
now formulate and prove the inductive step.

Suppose that for some $x_0\in \cS$, $\cS \setminus\{x_0\}=\cU_1
\cup \cU_2$, where $\cU_1$ and $\cU_2$ are disjoint and $\cU_1$ is
a finite union of finite trees $\cW_j$, $j=1,2,\dots, k$. Assume
that $Y_0 = x_0$ and $X_0 \in \cU_1$ (the argument is analogous
if the roles of $X$ and $Y$ are interchanged). Suppose without
loss of generality that $X_0 \in \cW_1$. Let $x_1\ne x_0$ be the
vertex of $\cW_1$ closest to $x_0$, and $\cW_1 \setminus \{x_1\} =
\cW_2 \cup \cW_3$, where $\cW_2$ and $\cW_3$ are disjoint, and
$\cW_3$ is the edge joining $x_0$ and $x_1$. We will first assume
that $\cW_2\ne \emptyset$. We will show that for some $p_3>0$ and
some stopping time $T_6< \infty$, with probability greater than
$p_3$, we either have $\d(X_{T_6}, Y_{T_6})\leq \eps$ or $Y_{T_6} =
x_1$ and $X_{T_6} \in \cW_2$.

If $\d(X_0, Y_0)\leq \eps$ then we can take $T_6 = 0$.

Assume that $\d(X_0, Y_0)> \eps$. Let $x_2\in \cW_3$ be the point
with $\d(x_0,x_2) = \eps/3$. Note that $\d(x_2,Y_0) \leq (1/2)
\d(x_2, X_0)$. Let $T_7 = \inf\{t\geq 0: X_t = x_2\}$ and $T_8 =
\inf\{t\geq 0: Y_t = x_2\}$. By Lemma 3.4, for some $s>0$,
 $$\P(T_8 < s) > \P(T_7 < s).$$
Hence, with probability $p_3>0$, $T_8< T_7$ and either $X$ and $Y$
have met by the time $T_8$, or $X$ is on the opposite side of
$Y_{T_8}$ in $\cS$ than $x_0$. Let $T_{9} = \inf\{t\geq T_8: Y_t =
x_1\}$ and $T_{10} = \inf\{t\geq T_8: Y_t = x_0\}$. Since $\cS$
contains only a finite number of finite trees $\cT_k$, there is an
upper bound on the edge length in any tree $\cT_k$, say,
$\rho<\infty$. This and the fact that $\d(x_2, x_0) =\eps/3$ imply
that $\P(T_9 < T_{10}) \geq p_4$ for some $p_4>0$ that may depend
on $\eps$. If the events $\{T_8< T_7\}$ and $\{T_9 < T_{10}\}$
hold then either $X$ and $Y$ have met by the time $T_9$ or
$Y_{T_9})= x_1$ and $X_{T_9})\in \cW_2$.

We note that if $\cW_2=\emptyset$ (i.e., $\cW_1 $ is a single
edge) then the same argument proves that for some $p_3>0$ and some
stopping time $T_6< \infty$, we have $\d(X_{T_6}, Y_{T_6})\leq \eps$
with probability greater than $p_3p_4$.

(iii) Let us rephrase the claim proved in step (ii). We have shown
that for some $p_5\df p_3p_4>0$ and some stopping time $T_6<
\infty$, with probability greater than $p_5$, we either have
$\d(X_{T_6}, Y_{T_6})\leq \eps$ or $Y_{T_6}$ and $X_{T_6}$ satisfy
the same assumptions as $Y_0$ and $X_0$, but relative to graphs
$\wt \cU_1 \df \cW_2$ and $\wt \cU_2 \df \cS\setminus (\{x_2\}
\cup \cW_2)$ in place of $\cU_1$ and $\cU_2$. Recall the claim
proved in part (i) of the proof and the final remark in step (ii).
Note that $\wt \cU_1$ has at least one edge less than $\cU_1$ so
by induction, we can repeat the inductive step (ii) a finite
number of times and show that with a probability $p_6>0$, $X$ and
$Y$ come within $\eps$ of each other before some time $t_1<
\infty$. It is easy to check that $p_6$ and $t_1$ can be chosen so
that they do not depend on the starting points of $X$ and $Y$. The
Markov property and induction can be used to show that $X$ and $Y$
have to come within $\eps$ of each other by the time $jt_1$ with
probability greater than $1-(1-p_6)^j$. We let $j\to\infty$ to see
that $X$ and $Y$ come within $\eps$ of each other at some finite
time with probability one. Since $\eps \in (0, r_0/3)$ is
arbitrary, the coupling is not shy.
 \qed

\bigskip
\noindent{\bf Example 3.6}. Suppose that $\cS$ is a tree with the
property that it has a ``backbone'' that is topologically a line,
with a finite or countable number of finite trees attached to it.
See Fig.~3.4 for an example.
Recall that we have assumed that each edge has length at least
$r_0>0$.

\bigskip \vbox{ \epsfxsize=3.0in
  \centerline{\epsffile{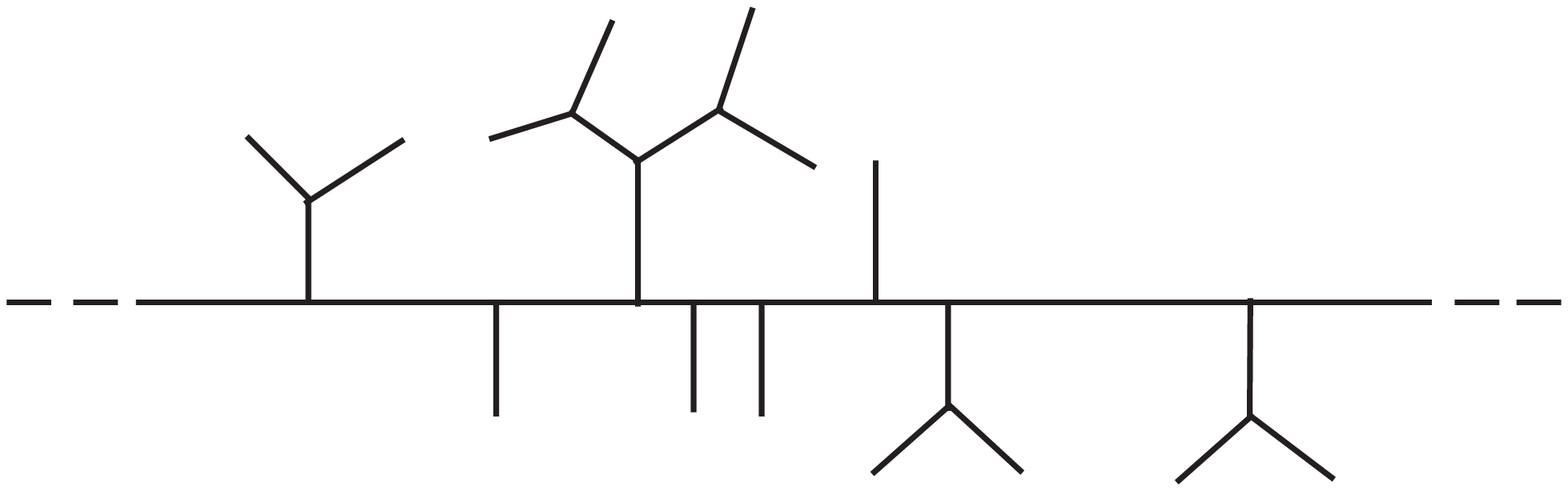}}

\centerline{Figure 3.4.} }
\bigskip

Let $\dots, x_{-2}, x_{-1}, x_0, x_1, x_2, \dots$ be the sequence
of points along the ``backbone'' $\cU$ where the side trees are
attached (the sequence can be finite or it can extend to infinity
in one or two directions). If the sequence extends to infinity in
both directions and
the graph is invariant under a non-constant shift,
 then, according to Example 3.3, there exists a
shy coupling.

Assume that
 \item{(i)} the diameters of the side trees are
uniformly bounded and
 \item{(ii)} $\{x_k\}$ does not extend to
infinity in both directions, or $\{x_k\}$ extends to infinity in
both directions but
the family $\{x_k\}$ is not shift-invariant, i.e., for every $c\ne
0$, there exists $x_j$ such that $x_j + c \ne x_k$ for all $k$.

We will show that under these assumptions there is no shy
coupling.

Let $Q$ be the ``projection'' of $\cS$ on $\cU$, i.e., $Q(x) = x$
for $x\in \cU$, and $Q(x) = x_k$ if $x$ belongs to a tree which is
attached to $\cU$ at $x_k$. We will identify the ``backbone''
$\cU$ with the real line so that we can think of $Q(X)$ and $Q(Y)$
as real-valued processes.
Let $X$ and $Y$ be a coupling of Brownian motions on $\cS$.
 If $X$ makes an excursion into a side tree then $Q(X)$ remains
constant on the excursion interval (including the endpoints).
Hence, $Q(X)$ and $Q(Y)$ are continuous processes. It is easy to
see that they are local martingales. Informally speaking, they are
Brownian motions frozen on some random intervals. The process $Z_t
= Q(X_t) - Q(Y_t)$ is also a continuous local martingale. Suppose
without loss of generality that $Z_0 >0$ and let
$$ T_0 =\inf\{t\geq 0: Z_t =0\} .
$$
Then $Z_{T_0\land t}$ is a non-negative local martingale and so it must
have an almost surely finite limit $Z_\infty$ on $\{T_0=\infty\}$.

We will show first that $\P(T_0=\infty)=\P(T_0=\infty \hbox{ and }
Z_\infty =0)$, in other words, $\P(T_0=\infty \hbox{ and }
Z_\infty >0)=0$. Let $B$ be a Brownian motion on $\R$ with
$B_0=0$. Consider a small $\eps>0$ and $t_0>0$, let $\delta \in
(0, \, \eps)$ be such that
$$
\P \left(\sup_{t\in[0,t_0]} |B_t| < \delta \right) < p_1/2,
\qquad \hbox{where } \ \
p_1:= \P \left(\sup_{t\in [0,t_0]}|B_{t} | <\eps \right).
$$

Let $T_1$ be the first time when all of the following conditions
hold: $\d(X_{T_1},\cU) \geq \eps$, the distance from $X_{T_1}$ to
any vertex of $\cS$ is greater than $\eps$, $Y_{T_1} \in \cU$, and
$\inf_k \d(Y_{T_1}, x_k) \geq \delta$ (the argument is analogous
if the roles of $X$ and $Y$ are interchanged). If $T_1<\infty$
then with probability $p_1/2$ or greater, $X$ will stay on the
same side tree over the interval $[T_1, T_1+t_0]$, while $Y$ will
move away from $Y_{T_1}$ by more than $\delta$ units over the same
time interval. Hence with probability $p_1/2$ or greater, $Z_t$
will have an oscillation of size at least $\delta$ over the
interval $[T_1, T_1+t_0]$. We proceed by induction. If $T_k <
\infty$ then we define $T_{k+1}=T_1 \circ \theta_{T_k+t_0} +
T_k+t_0$, where $\theta_\cdot$ is the usual Markovian shift
operator. Then with probability greater than $p_1/2$, $Z_t$ has an
oscillation of size at least $\delta$ over the interval $[T_{k+1},
T_{k+1}+t_0]$, independent of whether that happened over any
interval $[T_j, T_j+t_0]$, $j\leq k$. Hence, with probability one,
either $T_k= \infty$ for some $k$ or $Z_t$ has an infinite number
of oscillations of size $\delta$ over disjoint intervals of length
$t_0$, and, therefore in the latter case, $Z_t$ does not have a
limit as $t\to \infty$. Applying the above argument to a
decreasing sequence of $\{\eps_n, n\geq 1\}$ and a decreasing
sequence of $\{\delta_n, n \geq 1\}$ both tending to zero and
after deleting a null set from $\Omega$, we may and do assume that
for every $\omega \in \Omega$ and for every $\eps_n$, there is
some $N>1$
such that for every $j\geq N$, with $\eps_n$ and $\delta_j$ in place of
$\eps$ and $\delta$ above, either $T_k (\omega) =\infty$ for some
$k$ or $Z_t (\omega) $ does not have a limit as $t\to \infty$. The
processes $X$ and $Y$ are recurrent because the one-dimensional
Brownian motion is. Hence,
for every $x_j$,
 each one of them will enter the side tree attached to $\cU$ at
$x_j$ infinitely often. After deleting a null set from $\Omega$,
we may and do assume that the aforementioned property holds for
every $\omega\in \Omega$.

For $\omega \in \{T_0=\infty \hbox{ and } Z_\infty
>0 \}$, let
$$ c(\omega) =\lim_{t\to\infty} Z_t (\omega)>0 .
 $$
We choose an $x_j$, relative to $c(\omega)$, as follows. If $\{x_n
\}$ does not extend to $-\infty$ ($\infty$) then we let $x_j$ be
the leftmost (rightmost, resp.) point of the sequence. Otherwise
we fix an $x_j$ with the property that $x_j + c (\omega) \ne x_k$
for all $k$ (such an $x_j$ exists by assumption).
 Note that both $X(\omega)$ and $Y(\omega)$
enter the side tree attached to $\cU$ at $x_j$ infinitely often.
Hence one can find some $\eps>0$ from $\{\eps_n, n\geq 1\}$ and
$\delta>0$ from $\{\delta_n, n\geq 1\}$, and an increasing
sequence of random times $\{S_k, k\geq 1\}$ with $\lim_{k\to
\infty}S_k=\infty$ such that all of the following hold. One has
$\d(X_{S_k},\cU) \geq \eps$, the distance from $X_{S_k}$ to any
vertex of $\cS$ is greater than $\eps$, $Y_{S_k} \in \cU$, and
$\inf_n \d(Y_{S_k}, x_n) \geq \delta$ for every $k\geq 1$ (or the
statement will hold with the roles of $X$ and $Y$ interchanged).

Hence, all stopping times $\{T_k, k\geq 1\}$ defined in the
proceeding paragraph are finite. We have shown that this event
implies that $Z_\infty (\omega )$ does not exist. This
contradiction proves that $\P(T_0=\infty \hbox{ and }
Z_\infty>0)=0$ and therefore $\P(T_0=\infty)=\P(T_0=\infty \hbox{
and } Z_\infty =0)$.

Recall that we have assumed that all the edges have length at
least $r_0>0$. So on $\{T_0=\infty \hbox{ and } Z_\infty =0\}$, by
the recurrence of the one-dimensional Brownian motion, we have
$\liminf_{t\to \infty}\d(X_t, Y_t)=0$. We now only need to exam
$\omega \in \{T_0<\infty\}$ and to prove $X(\omega)$ and
$Y(\omega)$ will come arbitrarily close to each other.

Consider any $\eps\in (0, r_0/4)$, where $r_0>0$ is a lower bound
for the length of any edge in $\cS$. We want to show that with
probability one, there exists $t$ such that $\d(X_t, Y_t)\leq
\eps$. We have already proved that $\liminf_{t\to \infty} \d(X_t,
Y_t)=0$ on $\{T_0=\infty\}$. On $\{T_0<\infty\}$, at time $T_0$,
either $X_{T_0}=Y_{T_0}$ or one of processes $\{X_{T_0},
Y_{T_0}\}$ is at some $x_k$ and the other process is in a side
tree $\cT$ attached to $\cU$ at $x_k$. Without loss of generality,
assume that $X_{T_0}=x_k$ and $Y_{T_0}\in \cT$. If $\d(X_{T_0},
Y_{T_0})\leq\eps$ then we are done. Suppose that $\d(X_{T_0},
Y_{T_0})>\eps$. Let $z_0$ be the point at the edge $e$ of $\cT$
that is attached to $\cU$ with $\d(z_0, x_k)=\d(z_0, \cU)=\eps/4$.
Let $S_1$ be the first time after $T_0$ when $X_t = z_0$. By Lemma
3.4, with probability $p_2>0$, $X_t$ reaches $z_0$ after $T_0$
before $Y_t$ gets there. If this event occurs, both $X_{S_1}$ and
$Y_{S_1}$ will have distance at least $\eps/4$ away from $\cU$.
Let $R_1$ be the first time after $S_1$ when both $X_t$ and $Y_t$
are outside $\cT$. An argument analogous to that in parts (ii) and
(iii) of Example 3.5 shows that with probability $p_3>0$, the
processes $X$ and $Y$ will meet during the time interval $[S_1,
R_1]$. In other words, conditioning on $\{T_0<\infty\}$ and
$\d(X_{T_0}, Y_{T_0})>\eps$, with probability at least $p_4\df
p_2p_3>0$, the processes $X$ and $Y$ will meet between times $T_0$
and the first time $R_1$ when they are both outside $\cT$. We
define for $k\geq 2$,
$$ S_k=S_1\circ \theta_{R_{k-1}}+R_{k-1} \qquad \hbox{and} \qquad
R_k=R_1\circ \theta_{S_k}+S_k.
$$
By the strong Markov property of $(X, Y)$,
$$ \P\left( \inf_{t\in [0, R_k)} \d(X_t, \, Y_t)\geq \eps \right) \leq (1-p_4)^k.
$$
Letting $k\to\infty$, we get
$$ \P\left( \inf_{t\in [0, \infty)} \d(X_t, \, Y_t)
\geq \eps \right) \leq 0
$$
for every $\eps>0$ and thus $(X,Y)$ is not a shy coupling.
 \qed

\bigskip
\noindent{\bf Example 3.7}. Suppose that $\cS$ is composed of a
loop $\cU$ with a finite number of finite trees attached to it at
points $x_k$, and the family $\{x_k\}$ is not rotation invariant
in the following sense. We can assume without loss of generality
that $\cU$ is isometric to the unit circle. For every $c\ne 0$,
there exists $x_j$ such that $x_j e^{ic} \ne x_k$ for all $k$. See
Fig.~3.5 for an example. We will show that in this case there is
no shy coupling.

\bigskip \vbox{ \epsfxsize=1.5in
  \centerline{\epsffile{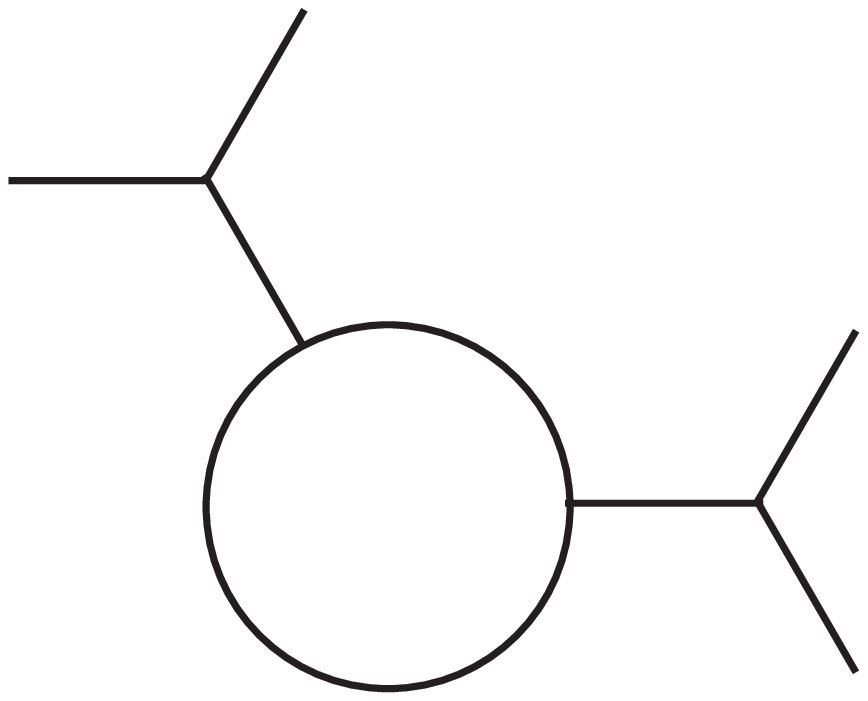}}

\centerline{Figure 3.5.} }
\bigskip

Our argument will be very similar to that in Example 3.6. Recall
the ``projection'' $Q$ from the previous example. We have $Q(x)=x$
for $x\in \cU$ and $Q(x) = x_k$ if $x$ belongs to a tree that is
attached to $\cU$ at $x_k$. Hence, $Q(X_t)$ may be regarded as a
continuous process on the unit circle. We now choose a (random)
continuous function $\Theta_X: [0,\infty) \to \R$ so that $Q(X_t)
= e^{i\Theta_X(t)}$ for all $t\geq 0$, in the complex notation. We
define $\Theta_Y$ in an analogous way. Note that $\Theta_X$ and
$\Theta_Y$ are martingales. Therefore, $Z_t \df \Theta_X(t) -
\Theta_Y(t)$ is also a martingale. We can now repeat the argument
from Example 3.6 to show that there does not exist a shy coupling.
 \qed

\bigskip
\noindent{\bf Example 3.8}. Examples 3.3 and 3.5-3.7 may appear to
suggest that if a graph has a vertex with degree 1 then a shy
coupling exists only if there exists an isometry of $I:\cS\to \cS$
with no fixed points. We will show that this is not the case. Our
example is illustrated in Fig.~3.6. In this case, every isometry
$I:\cS\to\cS$ has a fixed point. Nevertheless, we will show that
there is a shy coupling in $\cS$.

\bigskip \vbox{ \epsfxsize=3.0in
  \centerline{\epsffile{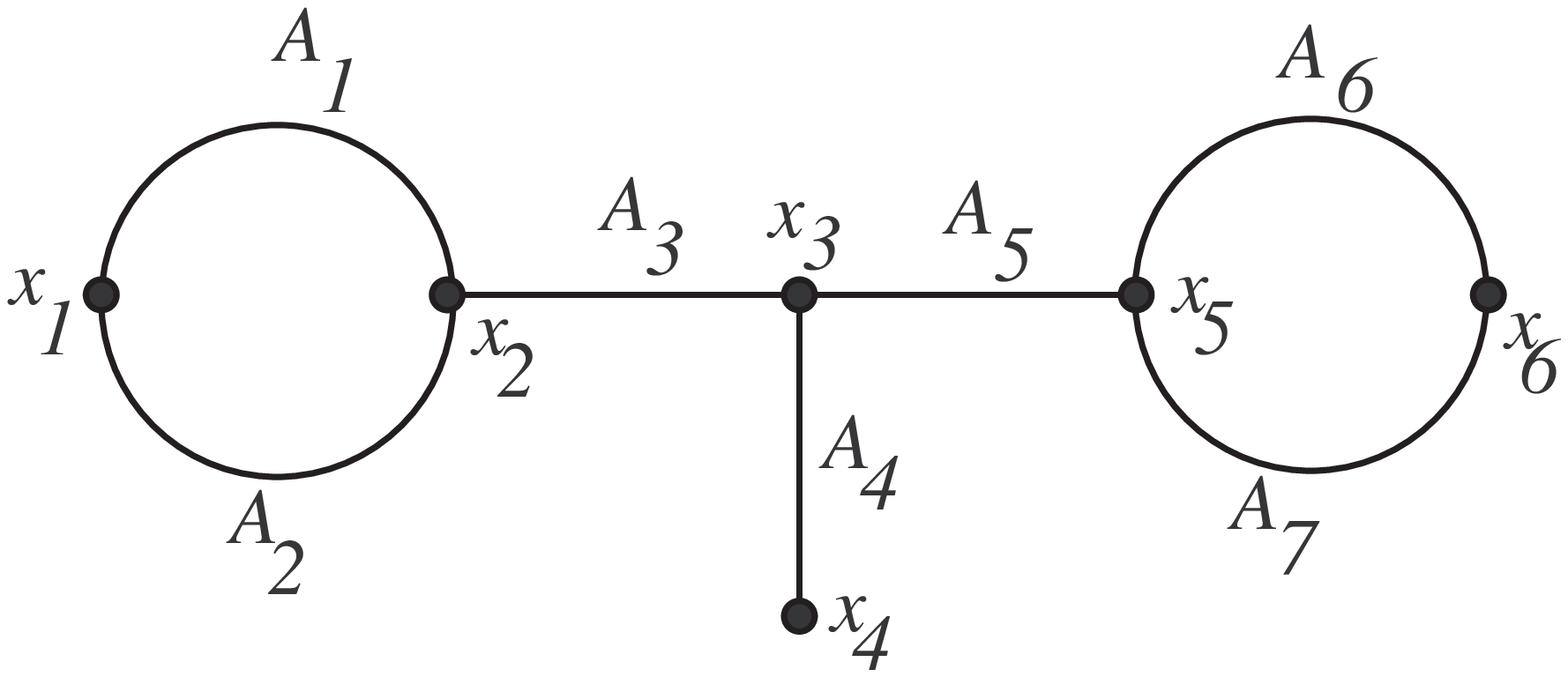}}

\centerline{Figure 3.6.} }
\bigskip

We will describe below the transition mechanism for $(X,Y)$ on
some random intervals of time. We will assume that the transition
probabilities of $(Y,X)$ are the same as those of $(X,Y)$. Hence,
there is no need to describe cases symmetric to those discussed
below, in the sense that the initial positions of $X$ and $Y$ are
interchanged.

Suppose that all edges $A_1, A_2, \dots, A_7$ have the same
length, say 1. We assume that $X_0 = x_2$ and $Y_0 = x_3$.

(i) Suppose that for some stopping time $T_1$ we have $X_{T_1} =
x_2$ and $Y_{T_1} = x_3$. Then we let $T_2 = \inf\{t\geq T_1: X_t
\notin (A_1 \cup A_2 \cup A_3)\setminus \{x_1\}\}$. Let $I: (A_1
\cup A_2 \cup A_3)\setminus \{x_1,x_3\} \to (A_3 \cup A_4 \cup
A_5)\setminus\{x_2,x_4,x_5\}$ be the one-to-one isometry
satisfying $I(x_2) = x_3$, $I(A_1) = A_3$, $I(A_2) = A_4$ and
$I(A_3) = A_5$. We let $Y_t = I(X_t)$ for $t\in[T_1,T_2]$. Note
that at the stopping time $T_2$, we have one of the following
configurations of the two particles: $(X_{T_2}, Y_{T_2})=(x_1,
x_4)$,  or $(X_{T_2}, Y_{T_2})= (x_1, x_2)$, or $(X_{T_2},
Y_{T_2}) = (x_3, x_5)$.

(ii) Suppose that for some stopping time $T_3$ we have $X_{T_3} =
x_1$ and $Y_{T_3} = x_4$. Let $\{X_t, t\in [T_3, T_4]\}$ be
Brownian motion on $\cS$ independent of the past with
$X_{T_3}=x_3$, where $T_4 = \inf\{t\geq T_3: X_t = x_2\}$. For
$t\in[T_3, T_4]$, we let $Y_t \in A_4$, with $\d(Y_t, x_4) =
\d(X_t, x_1)$. Note that $X_{T_4} = x_2$ and $Y_{T_4} = x_3$.

(iii) Suppose that for some stopping time $T_5$ we have $X_{T_5} =
x_1$ and $Y_{T_5} = x_2$. Let $\{Y_t, t\in [T_5, T_6]\}$ be
Brownian motion on $\cS$ independent of the past with
$Y_{T_5}=x_2$,
 where $T_6 = \inf\{t\geq T_5: Y_t = x_1 \hbox{ or } x_3\}$. We
label excursions of $Y$ from $x_2$ that stay in $A_3$ with marks
``1'' or ``2'', with equal probabilities, in such a way that the
label of any excursion is independent of all other labels. Then we
let $X_t$ be defined for $t\in [T_5, T_6]$ by $\d(X_t, x_1) =
\d(Y_t, x_2)$ and the following conditions. If $Y_t \in A_1$ then
$X_t \in A_2$, if $Y_t \in A_2$ then $X_t \in A_1$, if $Y_t \in
A_3$ and $t$ belongs to an excursion marked ``1'' then $X_t \in
A_1$, and if $Y_t \in A_3$ and $t$ belongs to an excursion marked
``2'' then $X_t \in A_2$. At time $T_6$ we have $(X_{T_6},
Y_{T_6})= (x_2,  x_1)$ or $(X_{T_6}, Y_{T_6}) =(x_2,  x_3)$.

Note that $\cS$ is symmetric with respect to the line containing
$A_4$. If for some stopping time $T_7$ we have $( X_{T_7},
Y_{T_7}) = (x_3, x_5)$, or $( X_{T_7}, Y_{T_7}) =(x_4,  x_6)$, or
$( X_{T_7}, Y_{T_7}) = (x_5,  x_6)$, or one of these conditions is
satisfied with the roles of $X$ and $Y$ interchanged, then we
define the coupling on an appropriate random interval in a way
analogous to that in (i)-(iii), using the symmetry of $\cS$.

The above definitions for the ``local'' behavior of the coupling
and the strong Markov property can now be used to define a process
$(X_t,Y_t)$ for all $t\geq 0$. It is easy to see that the stopping
times analogous to $T_1$, $T_3$ and $T_5$ will not have a finite
point of accumulation. It is also easy to check that almost surely
$ \d(X_t, Y_t) =1 $ for every $t>0$.
 \qed

\bigskip
\noindent {\bf 4. Reflected Brownian motion in Euclidean domains}.

This section is the closest in spirit to the papers and problems
which inspired the present research project. Suppose that
$D\subset \R^d$ is a bounded connected open set which is either
convex or has a $C^2$ boundary.
 We will consider couplings $(X,Y)$ of reflected Brownian motions
in $D$, defined as follows. Let $\n(x)$ denote the unit inward
normal vector at $x\in\prt D$. Let $B$ and $W$ be standard planar
Brownian motions with $B_0=W_0=0$ defined on the same probability
space and consider the following Skorohod equations,
 $$\eqalignno{
 X_t &= x_0 + B_t + \int_0^t  \n(X_s) dL^X_s, &(4.1)\cr
 Y_t &= y_0 + W_t + \int_0^t  \n(Y_s) dL^Y_s. &(4.2) }$$
Here $L^X$ is the local time of $X$ on $\prt D$, i.e., a
non-decreasing continuous process which does not increase when $X$
is in $D$: $\int_0^\infty \bone_{D}(X_t) dL^X_t = 0$, a.s.
Equation (4.1) has a unique pathwise solution $(X,L^X)$ such that
$X_t \in \ol D$ for all $t\geq 0$ (see [Ta] when $D$ is convex
domain and [LS] when $D$ is $C^2$ ). The ``reflected Brownian
motion'' $X$ is a strong Markov process. We point out that $B$ is
uniquely determined by $X$, and vice versa. The same remarks apply
to (4.2), so, as a pair, $(X, Y)$ is also strong Markov.

For a continuous semimartingale $M$, the symbol $\langle Z
\rangle$ will stand for its quadratic variation process. When
$M=(M^1, \cdots, M^d)$ and $N=(Z^1, \cdots , Z^d)$ are two
continuous $\R^d$-valued semimartingales we will use $\langle M,
N\rangle$ to denote $\sum_{i, j=1}^d \< M^i, N^j \>$. Note that
the matrix-valued process $( \< M^i, N^j \>)_{1\leq i, j\leq d}$
is non-negative definite and so $t\mapsto \< M , M\>_t$ is always
non-decreasing. For $a, b \in \R^d$, we use $a\cdot b$ to denote
the inner product between $a$ and $b$. We will use $\d(x, y)$ and
$|x-y|$ interchangeably for the Euclidean distance between $x, y
\in \R^d$.

\bigskip

\noindent{\bf Theorem 4.1}. {\sl Assume that $D\subset \R^d$ is a
bounded convex domain. Let $X$ and $Y$ be two reflecting Brownian
motion on $D$ given by (4.1)-(4.2).

\item{(i)} Suppose that there is a strictly increasing function
$\varphi$ with $\varphi (0)=0$ such that
$$ d\< |X-Y|^2\>_t \geq \varphi (|X_t-Y_t| )\, d t  \qquad \hbox{for }
t<\sigma_0,
$$
where $\sigma_0:=\inf \{t>0: X_t=Y_t\}$. Then $(X,Y)$ is not a shy
coupling.

\item{(ii)} Suppose
$D$ is strictly convex.
 Assume that $\< X-Y, X-Y\>_t$ (this is the same as $\< B-W,
B-W\>_t)$) has a sublinear growth rate as $t\to \infty$, that is,
$$\lim_{t\to \infty}  \< X-Y, \, X-Y\>_t/t=0 \qquad \hbox{almost surely}.
$$
Then $(X,Y)$ is not a shy coupling.

}

\bigskip
\noindent{\bf Proof}. Note that
$$ X_t-Y_t=X_0-Y_0+(B-W)+\int_0^t \n (X_s) dL^X_s -\int_0^t \n (Y_s) dL^Y_s
$$
is a semimartingale. Define $R_t:=|X_t-Y_t|^2$. By Ito's formula,
$$\eqalignno{
 d R_t &= 2 (X_t-Y_t) \cdot d(X_t-Y_t) + d\<X-Y, X-Y\>_t \cr
&= 2(X_t-Y_t) \cdot d(B_t-W_t) - 2 (Y_t-X_t)\cdot \n (X_t) dL^X_t
- 2(X_t-Y_t) \cdot \n (Y_t) dL^Y_t \cr & \hskip 0.2truein +
d\<X-Y, X-Y\>_t. &(4.3) \cr}
$$

\medskip

(i) Let
$a>0$ be a constant whose value will be
chosen in a moment and $f(r) = -r^{-a}$ for $r>0$. Then $f'(r) =
ar^{-a-1}>0$ and $f''(r) = (-a-1)a r^{-a-2}<0$ for $r>0$.
Define $U_t:= f(R_t)=f (|X_t-Y_t|^2)$. By Ito's formula, we have
$$
dU_t = f'(R_t) dR_t + {1\over 2} f'' (R_t) d\<R\>_t =dM_t+dV_t,
$$
where
$$ dM_t = 2 aR_t^{-a-1} 2(X_t-Y_t) \cdot d(B_t-W_t)
$$
and
$$\eqalignno{
dV_t &=  - 2 aR_t^{-a-1} \left( (Y_t-X_t) \cdot \n (X_t) dL^X_t
+ (X_t-Y_t)
\cdot
 \n (Y_t) dL^Y_t \right) \cr & \hskip 0.2truein  + a R_t^{-a-1}
d\<X-Y, X-Y\>_t -2 a(a+1) R_t^{-a-2} d\< |X-Y|^2\>_t  &(4.4) \cr}
$$
are the local martingale and bounded variation
parts,
 respectively. We claim that for every $\eps >0$,
$$ T_\eps := \inf \left\{ t>0:  \ |X_t-Y_t | \leq \eps \right\}
$$
is finite almost sure.
Suppose that $\P (T_\eps = \infty )>0$ for some $\eps >0$. We will
show that this leads to a contradiction.

Since $D$ is a convex domain, for $\sigma$-a.e. $x\in \partial D$,
$\n(x)$ is well defined and
 $$ (y-x) \cdot \n (x) \geq 0 \qquad \hbox{for every }
y
 \in \ol D.
$$
Note that the local time $L^X$ (respectively, $L^Y$) does not
increase when $X$ (respectively $Y$) is on a subset of $\partial
D$ having zero Lebesgue surface measure. Hence, (4.4) yields
$$
dV_t \leq   a R_t^{-a-1} d\<X-Y, X-Y\>_t -2 a(a+1) R_t^{-a-2} d\<
|X-Y|^2\>_t  . \eqno(4.5)
$$
Note that $ d\<X-Y, X-Y\>_t = d \<B-W, B-W\>_t \leq 4dt$. This,
(4.5) and the hypothesis in part (i) of this theorem imply that
on $\{T_\eps = \infty\}$,
$$\eqalignno{
dV_t &\leq   4 a R_t^{-a-1}dt  -2 a(a+1)
R_t^{-a-2} \varphi (\eps)   d t   \cr
&\leq  -2a R^{-a-2} \left( (a+1) \varphi (\eps )  - 2 R_t \right) dt \cr
& \leq -2a \eps^{-a-2}
\left( (a+1) \varphi (\eps )  -2 \hbox{\rm diam} (D) \right)dt . \cr}
$$
For a fixed $\eps>0$, we can find $a>0$ sufficiently large so that
for some $\lambda >0 $,
$$ dV_t \leq - \lambda  dt \qquad \hbox{for every } t>0
\hbox{ on } \{ T_\eps = \infty\}.
\eqno(4.6)
$$

The continuous local martingale $M$ is a time change of Brownian
motion. By the law of iterated logarithm for Brownian sample path,
for almost all $\omega \in \{T_\eps =\infty \}$, there is an
unbounded increasing sequence $\{t_k, k\geq 1\}$ such that
$\sup_{k\geq 1} | M_{t_k}(\omega) | <\infty $. This and (4.6)
imply that $U_{t_k} (\omega)= V_0(\omega)+M_{t_k}(\omega)+
V_{t_k}(\omega)$ tends to $-\infty$ as $k\to \infty$ on $\{T_\eps
=\infty \}$ a.s. Consequently, $|X_{t_k}-Y_{t_k}|$ goes to $0$ as
$k\to \infty$ on $ \{T_\eps =\infty \}$ a.s., which is a
contradiction. This proves that particles $X$ and $Y$ come
arbitrarily close to each other in finite time and, therefore, $X$
and $Y$ is not a shy coupling.

\medskip

(ii) Now assume
that $D$ is bounded and strictly convex
 and $\<X-Y, X-Y\>_t$ has a sublinear growth as $t\to \infty$.
The strict convexity implies (in fact, it is equivalent to) the
following condition. For every small $\eps >0$, there is a
constant $a_\eps>0$ such that
$$ (y-x) \cdot \n (x) \geq a_\eps \, |x-y| \qquad \hbox{for every }
x\in \partial D \hbox{ and } y \in \ol D \hbox{ with } |x-y| \geq
\eps. \eqno(4.7)
$$
Let $\sigma$ denote the surface measure on $\partial D$. Since
reflecting Brownian motion in $D$ is a recurrent Feller process,
it follows from the Ergodic Theorem that
$$ \lim_{t\to \infty} {L^X\over t}= {\sigma (\partial D ) \over 2 |D|}
=  \lim_{t\to \infty} {L^Y \over t} \qquad \hbox{almost surely.}
$$
For every $\eps >0$, define $T_\eps:=\inf\{t>0: |X_t-Y_t| \leq \eps\}$.
On $\{T_\eps =\infty\}$, we have from above
and (4.7)
that
$$\eqalignno{
& \liminf_{t\to \infty} {1\over t}
\left(-2 \int_0^t (Y_t-X_t)\cdot \n(X_t) dL^X_t
  - 2 \int_0^t (X_t-Y)\cdot \n(Y_t)  dL^Y_t + \<X-Y, X-Y\>_t \right) \cr
& \ \ \ \leq  \liminf_{t\to \infty} {1\over t} \left( -2 \eps
a_\eps L^X_t   - 2 \eps a_\eps  L^Y_t + \<X-Y, X-Y\>_t \right) =
-{ 2 \eps a_\eps \sigma ( \partial D) \over |D|}<0 . &(4.8)\cr}
$$
On the other hand, $M_t:=2\int_0^t (X_t-Y_t) \cdot d(B_t-W_t)$ is
a continuous martingale and thus is a time-change of
one-dimensional Brownian motion. By the law of iterated logarithm
for Brownian sample path, for almost all $\omega \in \{T_\eps
=\infty \}$, there is an unbounded increasing sequence $\{t_k,
k\geq 1\}$ such that $\sup_{k\geq 1} | M_{t_k}(\omega) | <\infty
$. This, (4.3) and (4.8) imply that $\lim_{t\to \infty}R_t=-\infty
$ a.s. on $\{T_\eps =\infty\}$. Since $R_t \geq 0$, we conclude
that $\P(T_\eps =\infty)=0$ for every $\eps>0$ and so $(X, Y)$ is
not a shy coupling.
 \qed

\bigskip
In the remainder of this section, we take $d=2$, but this is only
for notational convenience. We will show in the next example that
the method of proof of Theorem 4.1, based on the It\^o formula,
does not extend to arbitrary couplings. The example may have some
interest of its own. We will show in Theorem 4.3 below that, in
fact, there is no shy coupling of reflecting Brownian motions on
any bounded $C^1$-smooth strictly convex domain.

\bigskip

\noindent{\bf Example 4.2}. We will show that there exist planar
Brownian motions $B$ and $W$ with the property that $\d(B_t,W_t) =
\sqrt{2t + \d(B_0, W_0)^2}$ for $t\geq 0$, assuming that $B_0\ne
W_0$. In particular, the distance between the two processes grows
in a deterministic way.

Suppose that $(B, W)$ has the above mentioned property with $B_0$
and $W_0$ taking values in $\ol D$. Let $X$ and $Y$ be the
pathwise solutions
of (4.1)-(4.2) but with the above $B$ and $W$ in place of $x_0+B$
and $y_0+W$ there.
We have
$$ d\< |X_t-Y_t|^2\>= d\< |B_t-W_t|^2\>=0,
$$
while
$$ d\<X-Y, X-Y\>_t=d\<B-W, B-W\>_t =d \left( |B_t-W_t|^2\right)=2t.
$$
So $V_t$ in (4.4) becomes
$$\eqalignno{
dV_t &=  - 2 aR_t^{-a-1} \left( (Y_t-X_t) \cdot \n (X_t) dL^X_t
+ (X_t-Y_t)
\cdot
 \n (Y_t) dL^Y_t \right) \cr & \hskip 0.2truein  + 2 a R_t^{-a-1} d
t. \cr}
$$
Hence the method used in the proof of Theorem 4.1 does not work for
this coupling
 $(X , Y$). Moreover, since $|X_t-Y_t|$ grows deterministically
when both $X_t$ and $Y_t$ are away from the
boundary and decreases when one of them is on the boundary,
neither $|X_t-Y_t|$ nor any deterministic monotone
function of $|X_t-Y_t |$ is a submartingale or a supermartingale.

We now present the construction of $B$ and $W$ with the properties
mentioned above.
 Let $B$ be a Brownian motion in $\R^2$ starting from $x_0$. For a
vector $v=(a, b)\in \R^2$, we use $v^\perp$ to denote its
orthogonal vector $(b, -a)$. Let $y_0 \in \R^2$ be a point
different from $x_0$. Consider the following SDE for $W$ in $\R^2$
with $W_0=y_0$:
$$ dW_t= {1\over |W_t-B_t|^2} \left( \left( (W_t-B_t)\cdot dB_t \right) (W_t-B_t)
   - \left( (W_t-B_t)^\perp \cdot dB_t \right) (W_t-B_t)^\perp \right) .
$$
In words, at any given time $t>0$, $W_t$ takes
a synchronous step
 with $B_t$ along the direction $W_t-B_t$, while $W_t$ moves in the
opposite direction but with the same magnitude as $B_t$ along the
perpendicular direction $(W_t-B_t)^\perp$. The above SDE for $W$
has a unique solution up to $\tau:=\inf\{t>0: W_t=B_t\}$, since
the diffusion coefficients are $C^\infty$ up to that time. It can
be computed directly that $d\left( |W_t-B_t|^2 \right)= 2 dt$ and
consequently $|W_t-B_t|^2= |x_0-y_0|^2+2t$. So $\tau =\infty$.
It is standard to check
 that $W=(W^1, W^2)$ is a continuous local martingale with $\< W^i,
W^i\>_t=t$ for $i=1, 2$ and $\<W^1, W^2\>=0$. Therefore $W$ is a
Brownian motion in $\R^2$ starting from $y_0$.  \qed

\bigskip

We will show next that
in a $C^1$-smooth strictly convex domain $D$, every coupling of
reflecting Brownian motions on $\ol D$ must come arbitrarily close
to each other
 in finite time.

\bigskip

\noindent{\bf Theorem 4.3}. {\sl Suppose that $D$ is a bounded
convex planar domain with a
$C^1$-smooth
boundary that does not
contain any line segments. Then there does not exist a shy
coupling $(X,Y)$ of reflected Brownian motions in $D$.

}

\bigskip
\noindent{\bf Proof}. The idea of the proof is inspired by
differential games of pursuit (see [F]). We will show that with
positive probability, one of the particles will pursue the other
one in such a way that the distance between the two particles
decreases either because the diffusion component of the second
process does not move the second particle sufficiently fast or the
second particle hits the boundary and is pushed back towards the
first one.

\medskip

{\it Step 1}. We will define several constants $\eps_k$ in this
step. The definitions will be labeled (a), (b), (c), etc. Each of
these definitions is really a simple lemma asserting the existence
of a constant with stated properties. Since the proofs do not need
more than high school geometry, we omit most of the proofs. The
constants $\eps_k$ are defined relative to each other, but
$\eps_k$ may depend only on the values of $\eps_j$ for $j<k$.

(a) Let $\eps_0>0$ be so small that for every $x\in D$ with
$\d(x,\prt D) \leq \eps_0$ there exists a unique point in $\prt D$
whose distance from $x$ is minimal.

We make $\eps_0>0$ smaller, if necessary, so that the following is
true. Consider any point $y\in \prt D$ and let $CS_1$ be the
orthonormal coordinate system such that $y=0 \in \prt D$ and
$\n(0) $ lies on the second axis. Write $\n(x) = (\n_1(x),
\n_2(x))$. Then $|\n_1(x)| \leq \n_2(x)/100$ for $x \in \prt D
\cap \cB(0, \eps_0)$ in $CS_1$.

We fix an arbitrary $\eps_1 \in (0, \eps_0]$. It will suffice to
prove that for any $x_0,y_0 \in \ol D$, if $(X_0,Y_0)=(x_0,y_0)$
then, with probability one, there exists $t< \infty$ such that
$\d(X_t,Y_t) \leq \eps_1$.

(b) The angle between two vectors will be denoted $\angle
(\cdot\,,\,\cdot)$, with the convention that it takes values in
$(-\pi, \pi]$. Since $D$ is a bounded and strictly convex domain,
there exists $\eps_2 \in (0, \pi /2)$ such that for every $x\in
\prt D$ and $y\in \ol D$ satisfying $\d(x,y) \geq \eps_1/2$,
 $$ \angle (\n(x), y-x ) \in [ - \pi/2 +\eps_2, \, {\pi / 2}-\eps_2 ].
\eqno(4.9)
$$

(c) Let $\cL(x,r) $ be the cone spanned by $ \{\n (y)$, $y\in \prt
D\cap \cB(x, r)\}$. Since $\prt D$ is $C^1$-smooth, $\cL(x,r)$ is
a wedge. Hence, all linear combinations of vectors in $\cL(x,r)$
with non-negative coefficients belong to $\cL(x,r)$. An easy
approximation argument shows that if $X_t \in \cB(x,r)$ for all
$t\in(s,u)$ then $\int_s^u \n(X_t) dL^X_t \in \cL(x,r)$.

We will now choose $\eps_3\in(0,\eps_1/8)$. Consider vectors $\v$
and $\w$ satisfying the following conditions, relative to $x_0,
y_0$, and $\eps_3$.

If $\d(x_0,\prt D)\leq \eps_3$ then $ \v \in \cL(x_0, 2\eps_3)$
and $|\v| \leq \eps_3$. If $\d(x_0,\prt D)> \eps_3$ then $ \v =
0$.

If $\d(y_0,\prt D)\leq \eps_3$ then $ \w \in \cL(y_0, 2\eps_3)$
and $|\w| \leq \eps_3$. If $\d(y_0,\prt D)> \eps_3$ then $ \w =
0$.

We will show that (4.9) implies that we can find sufficiently
small $\eps_3>0$ so that the following is true. Suppose that $x_0,
y_0 \in \ol D$ with $\d(x_0, y_0) \geq \eps_1$. Assume that
$x_1\in \cB(x_0, 2\eps_3)$ and $y_1\in \cB(y_0, 2\eps_3)$. Then
 $$\d(x_1 + \v, y_1 +\w) \leq \d(x_1, y_1).\eqno(4.10)$$
To see this, choose some $x_2$ and $y_2$ so that the following
conditions hold.

If $\d(x_0,\prt D)\leq \eps_3$ then $x_2\in \prt D \cap \cB(x_0, 2\eps_3)$
with $ \v =c \n(x_2)$.
If $\d(x_0,\prt D)> \eps_3$ then $ x_2 = x_ 0$.

If $\d(y_0,\prt D)\leq \eps_3$ then $y_2\in \prt D \cap \cB(y_0, 2\eps_3)$
with $ \w =c \n(y_2)$.
If $\d(y_0,\prt D)> \eps_3$ then $ y_2 = y_ 0$.

Since $\eps_3< \eps_1 /8$,
$$ |x_2-y_2|\geq |x_0-y_0|-|x_0-x_2|-|y_0-y_2| \geq \eps_1 -4 \eps_3
    \geq \eps_1 /2.
$$
Thus by (4.9), we have
for sufficiently small $\eps_3>0$,
$$ \eqalignno{
&\d(x_1 + \v, y_1 +\w)^2 \cr &\leq |x_1-y_1|^2+|\v-\w|^2+2
(x_1-y_1)\cdot (\v-\w) \cr &\leq |x_1-y_1|^2+|\v-\w|^2+2
(x_2-y_2)\cdot (\v-\w) + 2(|x_1-x_2|+|y_1-y_2|)|\v-\w|\cr &\leq
|x_1-y_1|^2+ |\v -\w|^2 -2 (\sin \eps_2) |x_2-y_2| ( |\v|+|\w|) +
4 (8\eps_3)|\v-\w| \cr
 & \leq |x_1-y_1|^2 + (|\v|+|\w|) \left( 6\eps_3 -2 (\sin \eps_2)
(\eps_1/2)
 + 32 \eps_3 \right) \cr
 & =
 |x_1-y_1|^2 - (|\v|+|\w|) \left(
 \eps_1 \sin \eps_2
 -38 \eps_3 \right) &(4.11)\cr & \leq \d(x_1,
y_1)^2 .  \cr}
$$

(d) Since $D$ is bounded, we can find $\eps_4>0$ and $N<\infty$,
such that if $x_1, x_2, \dots$ is a sequence of points with
$x_1\in \ol D$, $\d(x_k,x_{k-1})\geq \eps_3/8$ and $|\angle(x_k -
x_{k-1}, x_{k+1} - x_k)| \leq \eps_4$ for all $k$ then $x_N\notin
\ol D$.

(e) We choose $\eps_5, \eps_6 >0$ so that the following is true.
Suppose that $x_0,y_0,x_1,y_1$ and $x_2$ satisfy the conditions
$\d(x_0,y_0)\geq \eps_1$, $\d(x_1,y_1)\geq \eps_1$ and
 $$\eqalignno{
 &|\angle (x_1- y_1, x_0-y_0)| \leq \eps_5,\cr
 &\d(x_1 ,
 x_0 + (\eps_3/4) (y_0-x_0)/\d(x_0,y_0) ) \leq \eps_6,\cr
 &\d(x_2 ,
 x_1 + (\eps_3/4) (y_1-x_1)/\d(x_1,y_1) ) \leq \eps_6.\cr}$$
Then $|\angle(x_1 - x_{0}, x_{2} - x_1)| \leq \eps_4$.

We make $\eps_6 $ smaller, if necessary, so that $\eps_6 <
\eps_3/8$.

(f) We make $\eps_6>0$ smaller, if necessary so that the following
holds. Suppose that $x_0,y_0,x_1$ and $y_1$ satisfy the conditions
$\d(x_0,y_0)\geq \eps_1$,
 $$\eqalignno{
 &\d(x_1 ,  x_0 + (\eps_3/4) (y_0-x_0)/\d(x_0,y_0) ) \leq
 \eps_6,&(4.12)\cr
 &\d(y_1 ,  y_0 + (\eps_3/4) (y_0-x_0)/\d(x_0,y_0) ) \leq \eps_6.&(4.13)\cr
}$$
Then $|\angle(x_1 - y_1, x_0 - y_0)| \leq \eps_5$.

(g) We can find $\eps_7,\eps_8 >0$ with the following properties.
Suppose that $x_0,y_0\in \ol D$, $x_1,x_2\in \R^2$,
$\d(x_0,y_0)\geq \eps_1$ and the following conditions are
satisfied.

If $\d(x_0,\prt D)\leq \eps_3$ then $ \v \in \cL(x_0, 2\eps_3)$
and $|\v| \leq \eps_3$. If $\d(x_0,\prt D)> \eps_3$ then $ \v =
0$.

If $\d(y_0,\prt D)\leq \eps_3$ then $ \w \in \cL(y_0, 2\eps_3)$
and $|\w| \leq \eps_3$. If $\d(y_0,\prt D)> \eps_3$ then $ \w =
0$.

Assume that
 $$\eqalignno{
 &\d(x_1 , x_0 + (\eps_3/4) (y_0-x_0)/\d(x_0,y_0) ) \leq
 \eps_8,&(4.14)\cr
 &\d(y_1 , y_0 + (\eps_3/4) (y_0-x_0)/\d(x_0,y_0) ) \geq \eps_6/2,&(4.15)\cr
 & \d(y_1, y_0) \leq \eps_3/4 +\eps_8. &(4.16) \cr
 }$$
Then $\d(x_1+\v, y_1+\w) \leq \d(x_0,y_0) -\eps_7$.

(h) It is easy to see from (4.11) that we can strengthen (4.10) as
follows. We can make $\eps_7>0$ smaller, if necessary, so that if
$|\v| \geq \eps_6/2$ or $|\w| \geq \eps_6/2$, and the assumptions
stated in Step 1(c) hold then
 $$\d(x_1 + \v, y_1 +\w) \leq \d(x_1, y_1)-2\eps_7.\eqno(4.17)$$

\medskip

\noindent{\it Step 2}. Suppose that $X_0, Y_0 \in \ol D$ with
$\d(X_0,Y_0) \geq \eps_1$. Consider the following events,
 $$\eqalign{
 F_1(t)&=\{
 \d(X_{t}, Y_{t}) \leq \d(X_0,Y_0) - \eps_7\},\cr
 F_2(t) & =
 \{|\angle (X_{t}- Y_{t}, X_0-Y_0)| \leq \eps_5\} ,\cr
 F_3(t) &=
 \{\d(X_t, X_0 + (\eps_3/4)
 (Y_0-X_0)/\d(X_0,Y_0) ) \leq \eps_6\land \eps_8\}, \cr
 F_4(t) & =\{ \d(X_{t}, Y_{t}) \leq \d(X_0,Y_0) +
 \eps_7/(4N)\},\cr
F_5(t) & = (F_1(t) \cup F_2(t) ) \cap F_3(t) \cap F_4(t).
 }$$
We will show in Step 4 that $\P( F_5(t_1)) > p_1$ for some $t_1,
p_1>0$ that do not depend on $X_0$ and $Y_0$.

Let $\eps_9 = \eps_7/(16N) \land \eps_3/8\land \eps_6/8\land
\eps_8/5$. Recall that $B$ and $W$ are Brownian motions with
$B_0=W_0=0$ driving $X$ and $Y$ in the sense of (4.1)-(4.2) and
let
 $$\eqalign{
 A_1(t)&=\left\{B_t\in \cB \left( (\eps_3/4) (Y_0-X_0)/\d(X_0,Y_0), \
 \eps_9 \right) \right\},\cr
 A_2(t)&=\left\{\sup_{ s\in[0,t]}
  | B_s - (s/t) B_t | \leq \eps_9 \right\},\cr
 A_3(t)&= \left\{ |W_t| \leq \eps_3/4 + \eps_9 \right\},\cr
 A_4(t) & = \left\{\sup_{ s\in[0,t]} | W_s - (s/t) W_t| \leq \eps_9
 \right\},\cr
 A_5(t) &= A_1(t) \cap A_2(t) \cap A_3(t) \cap A_4(t).
 }$$
We will argue in the rest of this step of the proof that
$\P(A_5(t_1) )> p_1$ for some $p_1,t_1>0$. In later steps, we will
show that $A_5(t)\subset F_5(t)$.

Recall that $B$ is a two-dimensional Brownian motion with $B_0=0$
and let $T_r = \inf\{t\geq 0: |B_t| > r\}$. Note that, by Brownian
scaling,
 $$ \P(T_r < t)
 = \P\left(\max_{0\leq s \leq t}|B_s| >r \right)
 = \P\left(  \max_{0\leq s \leq 1}|B_s| > r/\sqrt{t}  \right).
 $$
By the large deviations principle (see [RY], Ch. VIII, Thm. 2.11),
 $$\lim_{t/r^2 \to 0} (2t/r^2) \log \P(T_r < t\mid B_0=0)
 = -1.\eqno(4.18)$$

Let
 $$\eqalign{
 r_0 &= \eps_3/4, \cr
 A_6(t) &= \{ T_{r_0} < t\}, \cr
 A_7 & = \{B_{T_{r_0}} \in \cB ( (\eps_3/4) (Y_0-X_0)/\d(X_0,Y_0),
 \ \eps_9 /2)\}, \cr
 A_8(t) & = \left\{ \sup_{T_{r_0} \leq s \leq
 T_{r_0}+t} | B_s - B_{T_{r_0}} |  \leq \eps_9/2 \right\}.
 }$$
Clearly $A_6(t)$, $A_7$ and $A_8(t)$ are independent, and $ A_6(t)
\cap A_7 \cap A_8(t)\subset A_1(t)$. So
$$
\P(A_1(t)) \geq \P_6(A_6(t)) \cdot \P(A_7) \cdot \P(A_8(t)).
$$
 Note that
$\lim_{t\to 0} P(A_8(t)) = 1$ by the strong Markov
property applied at $T_{r_0}$ and the event $A_7$ is independent of $t$.
This together with (4.18) yields
 $$\liminf_{t \to 0}\,  t  \log \P(A_1(t))  \geq  -r_0^2/2= -(\eps_3/4)^2/2.$$
We obtain directly from (4.18) that
$$\limsup_{t \to 0} \, t \log \P(A^c_3(t)) \leq  -(\eps_3/4+\eps_9)^2/2.$$
This implies that for all sufficiently small $t>0$ we have
 $$\P(A_1(t) \cap A_3(t)) >0.\eqno(4.19)$$

The process $\{B_s -  (s/t)B_t, \, s\in[0,t]\}$ is Brownian bridge
with duration $t$ seconds, i.e., Brownian motion starting from 0
and conditioned to be at 0 at time $t$. If $\wt T_r$ denotes the
hitting time of $r$ by the absolute value of the Brownian bridge
then we have a formula analogous to (4.18), $\limsup_{t/r^2 \to 0}
(2t/r^2) \log \P(\wt T_r < t) \leq -1$. Thus
$$   \limsup_{t \to 0} t \log \P(A^c_2(t)) \leq  -\eps_9^2/2.
\eqno(4.20)
$$
Note that Brownian motion $B=\{B_t, t\geq 0\}$ is a Gaussian
process. Since for every $0\leq s<t$, $B_s-(s/ t) B_t$ and $B_t$
have zero covariance, the process $\{B_s -  (s/t)B_t, \,
s\in[0,t]\}$ is independent of $B_t$. Similar remarks apply to
$W$. Hence
$$
\P( A_1(t) \cap A_2(t))=  \P(A_1(t)) \P(A_2(t))
\quad \hbox{ and } \quad
\P( A_3(t) \cap A_4(t))=  \P(A_3(t)) \P(A_4(t)).
$$
This, (4.19) and (4.20)
 imply that
$$ \P(A_5(t_1) ) > p_1 \qquad \hbox{for some } t_1,p_1>0.
 \eqno (4.21)
$$
\medskip

\noindent{\it Step 3}. Let $R_t = \int_0^t \n(X_s) dL^X_s$. We
will show that if $A_1(t_0)\cap A_2(t_0)$ holds for some $t_0>0$
then $|R_s| \leq 4\eps_9$ for every $s\in[0,t_0]$.

Since $A_1(t_0)$ and $ A_2(t_0)$ hold, $|B_t| \leq \eps_3/4
+2\eps_9 \leq \eps_3/2$ for every $t\leq t_0$. Thus when $\d(X_0,
\prt D) \geq \eps_3$,  $X_0+B_s \in \ol D$ for all $s\in [0,t_0]$.
By the uniqueness of the solution to (4.1), $X_s =X_0+B_s$ for
$s\in [0,t_0]$, and $R_s = 0$ for all $s\in[0,t_0]$.

Suppose that $\d(X_0, \prt D) \leq \eps_3$. By the assumptions
made in Step 1(a), there exists a unique point $y\in \prt D$ with
the smallest distance to $X_0$. Let $CS_1$ be the orthonormal
coordinate system such that $y=0 \in \prt D$ and $\n(0) $ lies on
the second axis. Recall from Step 1(a) that $\n(x) = (\n_1(x),
\n_2(x))$ and $|\n_1(x)| \leq \n_2(x)/100$ for $x \in \prt D \cap
\cB(0, 3\eps_3)$ in $CS_1$. Write $R_t = (R^1_t, R^2_t)$. By the
opening remarks in Step 1(c), $R_t \in \cL(0,r)$ if $X_s\in
\cB(0,r)$ for all $s\leq t$. This implies that $|R^1_t| \leq
R^2_t/100$ if $X_s\in \cB(0,3\eps_3)$ for all $s\leq t$.

Let $T_1 = \inf\{t\geq 0: |R_{t}| > 4 \eps_9\}$. We will assume
that $T_1 < t_0$ and show that this leads to a contradiction.
Since $|B_t| \leq \eps_3/2$ for every $t\leq t_0$, we have $|B_t +
R_t| \leq \eps_3/2 + 4\eps_9 \leq \eps_3$ for every $t\leq T_1$.
Hence,
$$
|X_t| \leq | X_0 | + |B_t + R_t| \leq \eps_3 + \eps_3 =
2 \eps_3 \quad \hbox{ for  every } t\leq T_1.
$$
It follows that $R^2_t\geq 0$ for every $t\leq T_1$ and
$|R^1_{T_1}| \leq R^2_{T_1}/100$. Note that, by Step 1(a), the
slope of the tangent line at points in $\cB(0, 3\eps_3)\cap
\partial D$ is between $-1/100$ and $1/100$ and that $\d(X_{T_1} ,
X_0 + B_{T_1})=|R_{T_1}|=4\eps_9$. The last observation and the
fact that $X_{T_1} \in \prt D$ imply that $\d(X_0 + B_{T_1}, \prt
D)\geq (2/3) |R_{T_1}|=8\eps_9/3$. Since $A_1(t_0)$ and $
A_2(t_0)$ hold, $\d(X_0 +B_t, D) \leq 2\eps_9 $ for every $t\leq
t_0$ and in particular for $t=T_1$. This contradiction proves the
claim that $|R_s| \leq 4\eps_9$ for every $s\in[0,t_0]$.

\medskip

Let $\wt R_t = \int_0^t \n(Y_s) dL^Y_s$. We claim that if
$A_3(t_0)\cap A_4(t_0)$ holds for some $t_0>0$ then $|\wt R_s|
\leq \eps_3$ for every $s\in[0,t_0]$. To see this, observe that
$\d(Y_0 +W_t, D) \leq \eps_3/4 + 2\eps_9 \leq \eps_3/2$ for $t\leq
t_0$ and use that same argument as in the case of $R_t$.

\medskip

\noindent{\it Step 4}. Fix $t_0>0$. We will show that
$A_5(t_0)\subset F_5(t_0)$. Assume that $A_5(t_0)$ holds.

Since $A_1(t_0)$ holds, we have in view of Step 3,
 $$\eqalign{
 &\d(X_{t_0}, X_0 + (\eps_3/4)
 (Y_0-X_0)/\d(X_0,Y_0) ) \cr
&\leq \d(B_{t_0},  (\eps_3/4)
 (Y_0-X_0)/\d(X_0,Y_0) ) + |R_{t_0}|
 \leq \eps_9 + 4\eps_9
 \leq \eps_6\land \eps_8.}$$
In other words, $F_3(t_0)$ holds.

Since $A_1(t_0)$ and $A_3(t_0)$ hold we have, using simple
geometry,
 $$\d(X_0 + B_{t_0}, Y_0 + W_{t_0}) \leq \d(X_0,Y_0) + 3 \eps_9.\eqno(4.22)$$
We will apply (4.10) with $x_1 = X_0 + B_{t_0}$, $y_1 = Y_0 +
W_{t_0}$, $\v = R_{t_0}$ and $\w = \wt R_{t_0}$. We have assumed
that $A_5(t_0)$ holds so $|B_{t}| \leq \eps_3$ and $|W_{t}| \leq
\eps_3$ for $t\leq t_0$. This and Step 3 imply that for $t\leq
t_0$,
 $$\d(X_t, X_0) \leq |B_t| + |R_t| \leq 2 \eps_3,$$
and similarly $\d(Y_t, Y_0)  \leq 2 \eps_3$. By the opening
remarks in Step 1(c),
$$ \v = \int_0^{t_0} \n(X_t) dL^X_t \in \cL(X_0, 2 \eps_3)
\quad \hbox{ and } \quad \w = \int_0^{t_0} \n(Y_t) dL^Y_t
\in \cL(Y_0, 2 \eps_3).
$$
By Step 3, $|\v| \leq \eps_3$ and $|\w| \leq
\eps_3$. We have shown that all the conditions listed in Step 1(c)
are satisfied so we can apply (4.10) to obtain
 $$\d(X_{t_0} , Y_{t_0})  \leq \d(X_0,Y_0) + 3 \eps_9.$$
This proves that $F_4(t_0)$ holds.

It will now suffice to show that if $F_2(t_0)$ does not hold then
$F_1(t_0)$ does. Assume that $F_2(t_0)$ does not hold.

If all of the following conditions hold,
 $$\eqalignno{
  & | B_{t_0} -  (\eps_3/4)
 (Y_0-X_0 )/\d(X_0,Y_0)| \leq \eps_6/2\land c_8,&(4.23)\cr
 & | W_{t_0} -  (\eps_3/4) (Y_0-X_0)/\d(X_0,Y_0) |
 \leq \eps_6/2,&(4.24)\cr
 &\d(X_{t_0},  X_0 + B_{t_0} ) \leq \eps_6/2,&(4.25)\cr
 &\d(Y_{t_0}, Y_0 + W_{t_0}  ) \leq \eps_6/2,&(4.26)
 }$$
then (4.12) and (4.13) hold with $(x_1, y_1) = (X_{t_0}, Y_{t_0})$
and $(x_0, y_0)=(X_0, Y_0)$, and this implies $F_2(t_0)$, which is
a contradiction. Hence, at least one of the conditions
(4.23)-(4.26) must fail. The first of these conditions holds
because $A_1(t_0)$ is true. By Step 3, (4.25) holds.

Suppose that (4.26) fails. In view of (4.22), we can apply (4.17)
to $x_1 = X_0 + B_{t_0}$, $y_1 = Y_0 + W_{t_0}$, $\v = R_{t_0}$
and $\w = \wt R_{t_0}$ to obtain
 $$\d(X_{t_0} , Y_{t_0})  \leq \d(X_0,Y_0) + 3 \eps_9 -2\eps_7
 \leq \d(X_0,Y_0) - \eps_7.$$
Hence, we have $F_1(t_0)$ in this case.

Suppose that (4.24) fails. Then, in view of Step 3, (4.14)-(4.16)
hold with $(x_0, y_0)=(X_0, Y_0)$, $(x_1, y_1)=(X_0 + B_{t_0}, Y_0
+ W_{t_0})$, $\v = R_{t_0}$ and $\w = \wt R_{t_0}$ and we have
$$
\d(X_{t_0} , Y_{t_0})  \leq \d(X_0,Y_0) - \eps_7.
$$
Hence, $F_1(t_0)$ holds. This proves that $A_5(t_0)\subset F_5
(t_0)$.

\medskip

\noindent{\it Step 5}. Fix some $\eps_1\in (0,\eps_0)$ and let
$\eps_j$'s be defined relative to $\eps_1$ as in Step 1. Let
$\rho$ be the diameter of $D$ and let $N_0$ be an integer greater
than $4\rho/\eps_7$. Recall from (4.21) in Step 2 that for some
$p_1,t_1>0$ we have $\P(A_5(t_1) )> p_1$. Let
 $$S_1 = \inf\{t\geq 0: \d(X_t,Y_t) \leq \eps_1
 \hbox{  or  } F_5(t) \hbox{ holds}\}\land (2t_1).
 $$
By (4.21) and Step 4, $\P(S_1 \leq t_1) > p_1$. Recall that
$\theta$ stands for the usual Markov shift and define
$$ S_0=0 \qquad \hbox{and} \qquad
 S_k = S_1 \circ \theta_{S_{k-1}} + S_{k-1}\quad \hbox{for } k\geq 1.
$$
Recall integer $N$ defined in Step 1. By the strong Markov
property, with probability no less than $p_1^{2NN_0}>0$, we have
$S_k-S_{k-1} \leq  t_1$ for all $k \leq 2NN_0$.

We will argue that if $\bigcap_{k\leq 2NN_0} \{ S_k-S_{k-1} \leq
t_1\}$ holds then $\d(X_t, Y_t) \leq \eps_1$ for some $t \leq
2NN_0 t_1$. Assume otherwise. Then $F_5(S_1) \circ \theta
_{S_{k-1}}$ holds for every $k \leq 2NN_0$. In particular,
$F_4(S_1) \circ \theta _{S_{k-1}}$ holds for every $k \leq 2NN_0$.
Let $F_6(t) = F_2(t)\cap F_3(t)$. Since $F_5(t) \subset F_4(t)
\cap ( F_1(t) \cup (F_2(t)\cap F_3(t)))$, for every $k \leq
2NN_0$, at least one of the events $F_1(S_1) \circ \theta
_{S_{k-1}}$ and $F_6(S_1) \circ \theta _{S_{k-1}}$ holds.

Consider any $j\leq N_0$. If $F_6(S_1) \circ \theta _{S_{k-1}}$
holds for $k= 2jN, 2jN+1, \dots, 2(j+1)N -1$ then $X_{S_k}$'s and
$Y_{S_k}$'s satisfy the following conditions for $k= 2jN, 2jN+1,
\dots, 2(j+1)N -1$,
 $$\eqalign{
 \d\left (X_{S_{k+1}}, X_{S_k}+ (\eps_3/4)(Y_{S_k} - X_{S_k})/
\d (X_{S_k},Y_{S_k}) \right) &\leq \eps_6, \cr
 \left|\angle \left( X_{S_k} - Y_{S_k}, X_{S_{k+1}} - Y_{S_{k+1}}
 \right) \right|   &\leq \eps_5.
 }$$
This implies, by Step 1(e), that for $k= 2jN, 2jN+1, \dots,
2(j+1)N -2$,
 $$
 |\angle(X_{S_{k+1}} - X_{S_k}, X_{S_{k+2}} - X_{S_{k+1}})|
  \leq \eps_4.
 $$
Since $F_3(S_1) \circ \theta _{S_{k-1}}$ holds, we also have
$\d(X_{S_{k+1}}, X_{S_k}) \geq \eps_3/8$ for the same range of $k$ (to
see this, recall from Step 1(e) that $\eps_6 < \eps_3/8$). Hence
$X_{S_{2(j+1)N-2}}$ must be outside $\ol D$, according to the
definition of $\eps_4$ and $N$ in Step 1(d). Since $X$ always
stays inside $\ol D$, at least one of the events $F_6(S_1) \circ
\theta _{S_{k-1}}$ must fail for some $2jN \leq k \leq 2(j+1)N
-1$. Hence, at least one event $F_1(S_1) \circ \theta _{S_{k-1}}$
holds for some $2jN \leq k \leq 2(j+1)N -1$. Since $F_4(S_1) \circ
\theta _{S_{k-1}}$ holds for every $k \leq 2NN_0$, there is a
reduction of at least $\eps_7/2$ in the distance between $X$ and
$Y$ on every interval $[S_{2jN}, S_{2(j+1)N}]$, that is,
 $$\d \left( X_{S_{2jN}}, \ Y_{S_{2jN}} \right)
 \leq \d \left( X_{S_{2(j-1)N}}, \ Y_{S_{2(j-1)N}} \right)
 - {\eps_7\over 2} \qquad \hbox{for every } j\in \{ 1, \cdots,   N_0 \}.
 $$
Summing over $j$ we obtain
$$ \d \left( X_{S_{2N_0N}}, \ Y_{S_{2N_0N}} \right) \leq \d( X_0, \ Y_0 )
 - {N_0 \eps_7\over 2} \leq \d(X_0, Y_0)-2\rho <0.
$$
This contradiction proves our claim that
$$
\hbox{if  } \ \bigcap_{k\leq 2NN_0} \{ S_k-S_{k-1} \leq t_1\} \
\hbox{holds,  then } \ \d(X_t, Y_t) \leq \eps_1 \ \hbox{ for some
} \ t \leq 2t_1 NN_0 .
$$

We have shown that $\d(X_t,Y_t) \leq \eps_1$ for some $t\leq 2t_1
NN_0$ with probability greater than $p_2 \df p_1^{2NN_0}>0$. By
the Markov property, $\d(X_t,Y_t) \leq \eps_1$ for some $t\leq
2jt_1 NN_0$ with probability greater than $1- (1-p_2)^j$. To
complete the proof, it suffices to let $j\to \infty$.
 \qed

\bigskip

Two of the assumptions on the boundary of $D$ made in Theorem 4.3,
that it is convex with $C^1$-smooth boundary and it does not
contain any line segments, are convenient from the technical point
of view but most likely one can dispose of them with analysis more
refined than that in our proof.

\bigskip
\noindent{\bf Example 4.4}. Suppose that $D$ is the annulus
$\{x\in \R^2: 1<|x|<2\}$. The rotation of $D$ around $(0,0)$ with
an angle in $(0, 2\pi)$ is an isometry with no fixed points.
Hence, there exists a shy coupling of reflected Brownian motions
in this annulus (see Example 3.3).

\bigskip

There are many open problems concerning existence of shy couplings
but we find the following two questions especially intriguing.
Recall that $\cB(x,r)$ denotes the open ball with center $x$ and
radius $r$.

\bigskip
\noindent{\bf Open problems 4.5}. {\sl (i) Does there exist a shy
coupling of reflected Brownian motions in $\cB((0,0),3) \setminus
\cB((1,0),1)$?

(ii) Does there exist a shy coupling of reflected Brownian motions
in any simply connected planar domain?

}

\bigskip
We end this paper with a vague remark concerning a potential
relationship between shy couplings and an old and well known
problem of ``fixed points.'' Suppose that $\cS$ is a topological
space. If every continuous mapping $I: \cS \to \cS$ has a fixed
point, i.e., a point $x\in\cS$ such that $I(x) = x$, then we say
that $\cS$ has the fixed point property. One of the most famous
fixed point theorems is that of Brouwer---it asserts that a closed
ball in $\R^d$ has the fixed point property. Spheres obviously do
not have the fixed point property. Some of our results may suggest
that a shy coupling exists if and only if the state space does not
have the fixed point property. Example 3.8 applied to the graph
illustrated in Fig.~3.6 shows that this conjecture is false at
this level generality. It is possible, though, that a weaker form
of this assertion is true---we leave it as an open problem.

\vskip1truein

\centerline{REFERENCES}
\bigskip

\item{[BC1]} K.~Burdzy and Z.-Q.~Chen, Local time flow
related to skew Brownian motion. {\it Ann. Probab. \bf 29} (2001),
1693-1715.

\item{[BC2]} K.~Burdzy and Z.-Q.~Chen, Coalescence of
synchronous couplings. {\it Probab. Theory Rel. Fields \bf 123}
(2002), 553--578.

\item{[BCJ]} K.~Burdzy, Z.-Q.~Chen and P.~Jones, Synchronous
couplings of reflected Brownian motions in smooth domains.
Preprint, 2005.

\item{[BK]} K.~Burdzy and H.~Kaspi, Lenses in skew Brownian
flow. {\it Ann. Probab. \bf 32}  (2004), 3085--3115.

\item{[FW]} M.~Freidlin and A.~Wentzell, Diffusion processes
on graphs and the averaging principle. {\it Ann. Probab. \bf 21}
(1993), 2215-2245.

\item{[F]} A.~Friedman, {\it Differential Games}. Wiley, New
York, 1971.

\item{[HS]} J.M.~Harrison and L.A.~Shepp, On skew Brownian
motion. {\it Ann. Probab}. {\bf 9}  (1981), 309-313.

\item{[KS]} I.~Karatzas and S. E.~Shreve,  {\it Brownian Motion and
Stochastic Calculus}, Second edition. Springer, New York, 1991.

\item{[L]} T.~Lindvall, {\it Lectures on the Coupling
Method}. Wiley, New York,  1992.

\item{[LS]} P.~L. Lions and A.~S. Sznitman,  Stochastic
differential equations with reflecting boundary conditions. {\it
Comm. Pure Appl. Math.} {\bf 37}  (1984), 511-537.

\item{[RY]} D.~Revuz and M.~Yor, {\it Continuous Martingales
and Brownian Motion}. Springer, New York, 1991.

\item{[Ta]} H.  Tanaka, Stochastic differential equations with
reflecting boundary condition in convex regions. {\it Hiroshima
Math. J. \bf 9} (1979), 163-177.

\vskip1truein

\bigskip
\noindent I.B.: Department of Mathematics, Weizmann Institute of
Science, Rehovot 76100, Israel {\tt itai.benjamini@weizmann.ac.il}

\bigskip
\noindent K.B. and Z.C.: Department of Mathematics, Box 354350,
University of Washington, Seattle, WA 98115-4350, USA \hfill\break
{\tt burdzy@math.washington.edu, zchen@math.washington.edu}

\bye